\documentstyle[12pt]{article}

\catcode`\@=11

\ifcase\@ptsize
 \font\tenmsx=msxm10
 \font\sevenmsx=msxm7
 \font\fivemsx=msxm5
 \font\tenmsy=msym10
 \font\sevenmsy=msym7
 \font\fivemsy=msym5
\or
 \font\tenmsx=msxm10 scaled \magstephalf
 \font\sevenmsx=msxm8
 \font\fivemsx=msxm6
 \font\tenmsy=msym10 scaled \magstephalf
 \font\sevenmsy=msym8
 \font\fivemsy=msym6
\or
 \font\tenmsx=msxm10 scaled \magstep1
 \font\sevenmsx=msxm8
 \font\fivemsx=msxm6
 \font\tenmsy=msym10 scaled \magstep1
 \font\sevenmsy=msym8
 \font\fivemsy=msym6
\fi

\newfam\msxfam
\newfam\msyfam
\textfont\msxfam=\tenmsx  \scriptfont\msxfam=\sevenmsx
  \scriptscriptfont\msxfam=\fivemsx
\textfont\msyfam=\tenmsy  \scriptfont\msyfam=\sevenmsy
  \scriptscriptfont\msyfam=\fivemsy

\def\hexnumber@#1{\ifnum#1<10 \number#1\else
 \ifnum#1=10 A\else\ifnum#1=11 B\else\ifnum#1=12 C\else
 \ifnum#1=13 D\else\ifnum#1=14 E\else\ifnum#1=15 F\fi\fi\fi\fi\fi\fi\fi}

\def\msx@{\hexnumber@\msxfam}
\def\msy@{\hexnumber@\msyfam}
\mathchardef\boxdot="2\msx@00
\mathchardef\boxplus="2\msx@01
\mathchardef\boxtimes="2\msx@02
\mathchardef\square="0\msx@03
\mathchardef\blacksquare="0\msx@04
\mathchardef\centerdot="2\msx@05
\mathchardef\lozenge="0\msx@06
\mathchardef\blacklozenge="0\msx@07
\mathchardef\circlearrowright="3\msx@08
\mathchardef\circlearrowleft="3\msx@09
\mathchardef\rightleftharpoons="3\msx@0A
\mathchardef\leftrightharpoons="3\msx@0B
\mathchardef\boxminus="2\msx@0C
\mathchardef\Vdash="3\msx@0D
\mathchardef\Vvdash="3\msx@0E
\mathchardef\vDash="3\msx@0F
\mathchardef\twoheadrightarrow="3\msx@10
\mathchardef\twoheadleftarrow="3\msx@11
\mathchardef\leftleftarrows="3\msx@12
\mathchardef\rightrightarrows="3\msx@13
\mathchardef\upuparrows="3\msx@14
\mathchardef\downdownarrows="3\msx@15
\mathchardef\upharpoonright="3\msx@16

\mathchardef\downharpoonright="3\msx@17
\mathchardef\upharpoonleft="3\msx@18
\mathchardef\downharpoonleft="3\msx@19
\mathchardef\rightarrowtail="3\msx@1A
\mathchardef\leftarrowtail="3\msx@1B
\mathchardef\leftrightarrows="3\msx@1C
\mathchardef\rightleftarrows="3\msx@1D
\mathchardef\Lsh="3\msx@1E
\mathchardef\Rsh="3\msx@1F
\mathchardef\rightsquigarrow="3\msx@20
\mathchardef\leftrightsquigarrow="3\msx@21
\mathchardef\looparrowleft="3\msx@22
\mathchardef\looparrowright="3\msx@23
\mathchardef\circeq="3\msx@24
\mathchardef\succsim="3\msx@25
\mathchardef\gtrsim="3\msx@26
\mathchardef\gtrapprox="3\msx@27
\mathchardef\multimap="3\msx@28
\mathchardef\therefore="3\msx@29
\mathchardef\because="3\msx@2A
\mathchardef\doteqdot="3\msx@2B

\mathchardef\triangleq="3\msx@2C
\mathchardef\precsim="3\msx@2D
\mathchardef\lesssim="3\msx@2E
\mathchardef\lessapprox="3\msx@2F
\mathchardef\eqslantless="3\msx@30
\mathchardef\eqslantgtr="3\msx@31
\mathchardef\curlyeqprec="3\msx@32
\mathchardef\curlyeqsucc="3\msx@33
\mathchardef\preccurlyeq="3\msx@34
\mathchardef\leqq="3\msx@35
\mathchardef\leqslant="3\msx@36
\mathchardef\lessgtr="3\msx@37
\mathchardef\backprime="0\msx@38
\mathchardef\risingdotseq="3\msx@3A
\mathchardef\fallingdotseq="3\msx@3B
\mathchardef\succcurlyeq="3\msx@3C
\mathchardef\geqq="3\msx@3D
\mathchardef\geqslant="3\msx@3E
\mathchardef\gtrless="3\msx@3F
\mathchardef\sqsubset="3\msx@40
\mathchardef\sqsupset="3\msx@41
\mathchardef\vartriangleright="3\msx@42
\mathchardef\vartriangleleft="3\msx@43
\mathchardef\trianglerighteq="3\msx@44
\mathchardef\trianglelefteq="3\msx@45
\mathchardef\bigstar="0\msx@46
\mathchardef\between="3\msx@47
\mathchardef\blacktriangledown="0\msx@48
\mathchardef\blacktriangleright="3\msx@49
\mathchardef\blacktriangleleft="3\msx@4A
\mathchardef\vartriangle="3\msx@4D
\mathchardef\blacktriangle="0\msx@4E
\mathchardef\triangledown="0\msx@4F
\mathchardef\eqcirc="3\msx@50
\mathchardef\lesseqgtr="3\msx@51
\mathchardef\gtreqless="3\msx@52
\mathchardef\lesseqqgtr="3\msx@53
\mathchardef\gtreqqless="3\msx@54
\mathchardef\Rrightarrow="3\msx@56
\mathchardef\Lleftarrow="3\msx@57
\mathchardef\veebar="2\msx@59
\mathchardef\barwedge="2\msx@5A
\mathchardef\doublebarwedge="2\msx@5B
\mathchardef\angle="0\msx@5C
\mathchardef\measuredangle="0\msx@5D
\mathchardef\sphericalangle="0\msx@5E
\mathchardef\varpropto="3\msx@5F
\mathchardef\smallsmile="3\msx@60
\mathchardef\smallfrown="3\msx@61
\mathchardef\Subset="3\msx@62
\mathchardef\Supset="3\msx@63
\mathchardef\Cup="2\msx@64

\mathchardef\Cap="2\msx@65

\mathchardef\curlywedge="2\msx@66
\mathchardef\curlyvee="2\msx@67
\mathchardef\leftthreetimes="2\msx@68
\mathchardef\rightthreetimes="2\msx@69
\mathchardef\subseteqq="3\msx@6A
\mathchardef\supseteqq="3\msx@6B
\mathchardef\bumpeq="3\msx@6C
\mathchardef\Bumpeq="3\msx@6D
\mathchardef\lll="3\msx@6E

\mathchardef\ggg="3\msx@6F

\mathchardef\circledS="0\msx@73
\mathchardef\pitchfork="3\msx@74
\mathchardef\dotplus="2\msx@75
\mathchardef\backsim="3\msx@76
\mathchardef\backsimeq="3\msx@77
\mathchardef\complement="0\msx@7B
\mathchardef\intercal="2\msx@7C
\mathchardef\circledcirc="2\msx@7D
\mathchardef\circledast="2\msx@7E
\mathchardef\circleddash="2\msx@7F
\def\ulcorner{\delimiter"4\msx@70\msx@70 }
\def\urcorner{\delimiter"5\msx@71\msx@71 }
\def\llcorner{\delimiter"4\msx@78\msx@78 }
\def\lrcorner{\delimiter"5\msx@79\msx@79 }
\def\yen{\mathhexbox\msx@55 }
\def\checkmark{\mathhexbox\msx@58 }
\def\circledR{\mathhexbox\msx@72 }
\def\maltese{\mathhexbox\msx@7A }
\mathchardef\lvertneqq="3\msy@00
\mathchardef\gvertneqq="3\msy@01
\mathchardef\nleq="3\msy@02
\mathchardef\ngeq="3\msy@03
\mathchardef\nless="3\msy@04
\mathchardef\ngtr="3\msy@05
\mathchardef\nprec="3\msy@06
\mathchardef\nsucc="3\msy@07
\mathchardef\lneqq="3\msy@08
\mathchardef\gneqq="3\msy@09
\mathchardef\nleqslant="3\msy@0A
\mathchardef\ngeqslant="3\msy@0B
\mathchardef\lneq="3\msy@0C
\mathchardef\gneq="3\msy@0D
\mathchardef\npreceq="3\msy@0E
\mathchardef\nsucceq="3\msy@0F
\mathchardef\precnsim="3\msy@10
\mathchardef\succnsim="3\msy@11
\mathchardef\lnsim="3\msy@12
\mathchardef\gnsim="3\msy@13
\mathchardef\nleqq="3\msy@14
\mathchardef\ngeqq="3\msy@15
\mathchardef\precneqq="3\msy@16
\mathchardef\succneqq="3\msy@17
\mathchardef\precnapprox="3\msy@18
\mathchardef\succnapprox="3\msy@19
\mathchardef\lnapprox="3\msy@1A
\mathchardef\gnapprox="3\msy@1B
\mathchardef\nsim="3\msy@1C
\mathchardef\napprox="3\msy@1D
\mathchardef\varsubsetneq="3\msy@20
\mathchardef\varsupsetneq="3\msy@21
\mathchardef\nsubseteqq="3\msy@22
\mathchardef\nsupseteqq="3\msy@23
\mathchardef\subsetneqq="3\msy@24
\mathchardef\supsetneqq="3\msy@25
\mathchardef\varsubsetneqq="3\msy@26
\mathchardef\varsupsetneqq="3\msy@27
\mathchardef\subsetneq="3\msy@28
\mathchardef\supsetneq="3\msy@29
\mathchardef\nsubseteq="3\msy@2A
\mathchardef\nsupseteq="3\msy@2B
\mathchardef\nparallel="3\msy@2C
\mathchardef\nmid="3\msy@2D
\mathchardef\nshortmid="3\msy@2E
\mathchardef\nshortparallel="3\msy@2F
\mathchardef\nvdash="3\msy@30
\mathchardef\nVdash="3\msy@31
\mathchardef\nvDash="3\msy@32
\mathchardef\nVDash="3\msy@33
\mathchardef\ntrianglerighteq="3\msy@34
\mathchardef\ntrianglelefteq="3\msy@35
\mathchardef\ntriangleleft="3\msy@36
\mathchardef\ntriangleright="3\msy@37
\mathchardef\nleftarrow="3\msy@38
\mathchardef\nrightarrow="3\msy@39
\mathchardef\nLeftarrow="3\msy@3A
\mathchardef\nRightarrow="3\msy@3B
\mathchardef\nLeftrightarrow="3\msy@3C
\mathchardef\nleftrightarrow="3\msy@3D
\mathchardef\divideontimes="2\msy@3E
\mathchardef\varnothing="0\msy@3F
\mathchardef\nexists="0\msy@40
\mathchardef\mho="0\msy@66
\mathchardef\thorn="0\msy@67
\mathchardef\beth="0\msy@69
\mathchardef\gimel="0\msy@6A
\mathchardef\daleth="0\msy@6B
\mathchardef\lessdot="3\msy@6C
\mathchardef\gtrdot="3\msy@6D
\mathchardef\ltimes="2\msy@6E
\mathchardef\rtimes="2\msy@6F
\mathchardef\shortmid="3\msy@70
\mathchardef\shortparallel="3\msy@71
\mathchardef\smallsetminus="2\msy@72
\mathchardef\thicksim="3\msy@73
\mathchardef\thickapprox="3\msy@74
\mathchardef\approxeq="3\msy@75
\mathchardef\succapprox="3\msy@76
\mathchardef\precapprox="3\msy@77
\mathchardef\curvearrowleft="3\msy@78
\mathchardef\curvearrowright="3\msy@79
\mathchardef\digamma="0\msy@7A
\mathchardef\varkappa="0\msy@7B
\mathchardef\hslash="0\msy@7D
\mathchardef\hbar="0\msy@7E
\mathchardef\backepsilon="3\msy@7F
\def\Bbb{\ifmmode\let\next\Bbb@\else
 \def\next{\errmessage{Use \string\Bbb\space only in math mode}}\fi\next}
\def\Bbb@#1{{\Bbb@@{#1}}}
\def\Bbb@@#1{\fam\msyfam#1}

\font\teneuf=eufm10
\font\seveneuf=eufm7
\font\fiveeuf=eufm5
\newfam\euffam
\textfont\euffam=\teneuf
\scriptfont\euffam=\seveneuf
\scriptscriptfont\euffam=\fiveeuf
\def\frak{\relaxnext@\ifmmode\let\next\frak@\else
 \def\next{\Err@{Use \string\frak\space only in math mode}}\fi\next}
\def\goth{\relaxnext@\ifmmode\let\next\frak@\else
 \def\next{\Err@{Use \string\goth\space only in math mode}}\fi\next}
\def\frak@#1{{\frak@@{#1}}}
\def\frak@@#1{\noaccents@\fam\euffam#1}

\def\frak{\ifmmode\let\next\frak@\else
 \def\next{\errmessage{Use \string\frak\space only in math mode}}\fi\next}
\def\frak@#1{{\frak@@{#1}}}
\def\frak@@#1{\fam\euffam#1}

\newcommand{\nc}{\newcommand}

\newcommand{\colvec}[2]{\left  ( \begin{array}{cc} #1  \\
     #2  \end{array} \right ) }

\newcommand{\bla}{\phantom{bbbbb}}
\newcommand{\onebl}{\phantom{a} }
\newcommand{\eqdef}{\;\: {\stackrel{ {\rm def} }{=} } \;\:}

\newcommand{\half}{ {\frac{1}{2} } }
\newcommand{\vol}{ \,{\rm vol}\, }

\newcommand{\beq}{\begin{equation}}
\newcommand{\eeq}{\end{equation}}
\newcommand{\beqst}{\begin{equation*}}
\newcommand{\eeqst}{\end{equation*}}
\newcommand{\barr}{\begin{array}}
\newcommand{\earr}{\end{array}}
\newcommand{\beqar}{\begin{eqnarray}}
\newcommand{\eeqar}{\end{eqnarray}}
\newtheorem{theorem}{Theorem}[section]
\newtheorem{corollary}[theorem]{Corollary}
\newtheorem{lemma}[theorem]{Lemma}
\newtheorem{prop}[theorem]{Proposition}
\newtheorem{definition}[theorem]{Definition}
\newtheorem{remit}[theorem]{Remark}

\newcommand{\RR}{{\Bbb R }}
\newcommand{\CC}{{\Bbb C }}

\newcommand{\PP}{ {\Bbb P } }

\newcommand{\cala}{{\mbox{$\cal A$}}}
\newcommand{\calb}{{\mbox{$\cal B$}}}

\newcommand{\cald}{{\mbox{$\cal D$}}}

\newcommand{\calf}{{\mbox{$\cal F$}}}

\newcommand{\cali}{{\mbox{$\cal I$}}}
\newcommand{\calj}{{\mbox{$\cal J$}}}

\newcommand{\calm}{{\mbox{$\cal M$}}}

\newcommand{\calo}{{\mbox{$\cal O$}}}

\newcommand{\calr}{{\mbox{$\cal R$}}}
\newcommand{\cals}{{\mbox{$\cal S$}}}

\def\a{\alpha}
\def\b{\beta}
\def\g{\gamma}
\def\d{\delta}
\def\e{\epsilon}
\def\z{\zeta}

\def\t{\theta}
\def\k{\kappa}
\def\l{\lambda}
\def\m{\mu}
\def\n{\nu}

\def\s{\sigma}

\def\G{\Gamma}

\def\T{\Theta}
\def\L{\Lambda}

\def\O{\Omega}

\newcommand{\renorm}{{ \setcounter{equation}{0} }}

\nc{\Imm}{ {\rm Im} }
\nc{\Si}{\Sigma}
\nc{\si}{\sigma}

\nc{\liek}{  \mbox{\bf k}  }
\nc{\liet}{  \mbox{\bf t}  }

\nc{\lietp}{ {\liet}^\perp  }
\nc{\liets}{ {\liet}^* }
\nc{\lieks}{ {\liek}^* }
\nc{\hk}{H^*_K}
\nc{\hht}{H^*_T}
\nc{\pisk}{\Pi_*}
\nc{\pist}{\Pi_*}
\nc{\pis}{\Pi_*}

\nc{\be}{\beta}

\nc {\Om}{\Omega}
\nc{\om}{\omega}

\nc{\diag}{ {\rm diag} }

\nc{\lrar}{\longrightarrow}
\nc{\Proof}{ \noindent{\em Proof:} }
\nc{\Cok}{ {\rm Cok} }

\nc{\ie}{\cali^\e}
\nc{\xred}{\calm_X }

\nc{\efo}{e_F^{(0)} }
\nc{\inpr}[1]{ \langle #1 \rangle }
\nc{\zloc}{{ \mu^{-1}(0) } }
\nc{\fk}{F_K}
\nc{\ft}{F_T}

\nc{\eva}[2]{ { #1(#2)   } }
\nc{\evab}[2]{ { #2(#1)   } }

\nc{\epin}[1]{e^{ i \inpr{ #1} }  }
\nc{\emin}[1]{e^{ -i \inpr{ #1} }  }
\nc{\epinev}[2]{e^{ i \eva{ #1}{#2  } }  }
\nc{\eminev}[2]{e^{ -i \eva{ #1}{#2  } }  }
\nc{\epinevb}[2]{e^{ i \evab{ #1}{#2  } }  }
\nc{\eminevb}[2]{e^{ -i \evab{ #1}{#2  } }  }

\nc{\omr}{\om_0}
\nc{\intk}{\int_{\phi \in \liek} }
\nc{\intt}{\int_{\psi \in \liet} }
\nc{\tintk}{\int_{z \in \lieks} }
\nc{\tintt}{\int_{y \in \liets} }
\nc{\bo}{B_0}
\nc{\hl}{H_\l}
\nc{\vl}{V_\l}
\nc{\mub}{\mu_{\hat{\beta} } }
\nc{\hatb}{ { \hat{\beta} } }
\nc{\ftpis} {\ft (\pist \si) }
\nc{\thb}{T_\hatb}
\nc{\stg}{ e^{- |y|^2/{2 \e} }  }

\nc{\gse}{g_{\e} }
\nc{\gsoe}{g_{\e^{-1} } }
\nc{\zeb}{{\z'}_\b}
\nc{\bom}{ {\bar{\om} } }
\nc{\va}{V^j}
\nc{\eeth   }{e^{\e \T} e^{i\om_0} }
\nc{\gtsoe}{{\widetilde{\gsoe} } }
\nc{\wn}{{ |W|}}
\nc{\bfj} {\b_{F, j} }
\nc{\bfjw}{ \b_{F, j}^{\L} }
\nc{\fj}{ {F,j} }
\nc{\hh}{H}
\nc{\hhh}{H}
\nc{\ggh}{G}
\nc{\gggh}{G}
\nc{\hmo}{{H^0_\m} }
\nc{\cm}{ {\rm Crit} (\mu_T) }
\nc{\tb}{ {\tilde{B} } }
\nc{\tbg}{\tilde{B}^{(\g)} }
\nc{\tbo}{\tb_0}
\nc{\qeta}{ {Q^\eta} }
\nc{\hatf}{ \hat{f} }
\nc{\hatg}{ \hat{g} }
\nc{\fhtpl}{\frac{1} {(2 \pi)^{l/2} } }
\nc{\fhtps}{\frac{1} {(2 \pi)^{s/2} } }
\nc{\htpl}{ {(2 \pi)^{l/2} } }
\nc{\htps}{ {(2 \pi)^{s/2} } }
\nc{\ftpl}{\frac{1} {(2 \pi)^{l} } }
\nc{\ftps}{\frac{1} {(2 \pi)^{s} } }
\nc{\tpl}{ {(2 \pi)^{l} } }
\nc{\tps}{ {(2 \pi)^{s} } }
\nc{\dbyd}[1]{\frac{\partial}{\partial #1} }
\nc{\dbydl}[1]{{\partial}/{\partial #1} }
\nc{\dell}{\bigtriangleup}
\nc{\indd}{n}
\nc{\Indd}{m}
\nc{\gnx}{\xi}
\nc{\gna}{\a}
\nc{\gnp}{P}
\nc{\anr}{a_n^{(r)} }
\nc{\caljab}{ {\calj}^{\rm ab} }
\nc{\tpsi}{\tilde{\psi} }
\nc{\ttfa}{ { \tilde{\tau}_{F,\a} } }
\nc{\tga}{\tilde{\Gamma}}
\nc{\barl}{ {\bar{\l} } }

\nc{\tfaj}{\tau_{F, \a,J}^\eta  }
\nc{\reso}{  {  {\rm Res}^{\L,\testf}    }  }
\nc{\gdist}[1]{\frac
{e^{i #1 (\psi)}   } {\prod_{j = 1}^N \beta_j(\psi) } [d \psi] }
\nc{\gdistl}[1]{
{e^{i #1 (\psi)} [d \psi]  }/ {\prod_{j = 1}^N \beta_j(\psi) }  }
\nc{\limpl}{ {  \lim_{s \to 0^+} }}
\nc{\barb}{{\bar{\beta}} }
\nc{\hbp}{H'_\barb}
\nc{\veee}{v}
\nc{\vt}{\xi}
\nc{\treso}{ { {\rm Res} } }
\nc{\cmon}{\CC^{{\rm mon}, \xi} (\liet) }
\nc{\cmons}{\CC^{{\rm mon}, \xi}_s (\liet) }
\nc{\mar}{ }
\nc{\iins}{ {i} }
\nc{\nusym}{\varpi}
\nc{\rf}{r_F}
\nc{\rfe}{r_F^\eta}

\nc{\htestf}{\hat{\chi}}
\nc{\htestfe}{\hat{\chi}_\e}
\nc{\testf}{{\chi}}
\nc{\testfe}{{\chi}_\e}
\nc{\tx}{\tau_\xi}
\nc{\reslch}{{\rm Res}^{\L,\testf} \,}
\nc{\calje}{ { {\calj}^\e} }
\nc{\zray}{\rho}
\nc{\resrlch}{{\rm Res}^{\zray, \L,\testf} \, }
\nc{\rbare}{r^{{\eta} } }
\nc{\lasub}{\L}
\nc{\ress}{{\rm Res} \,}
\nc{\ck}{C_K}
\nc{\je}{ { {\calj}^\e} }
\nc{\deleps}{\bigtriangleup_\e}
\nc{\bareta}{\bar{\eta} }

\nc{\gu}{\G_u}
 \nc{\guo}{\G_u^0}
\nc{\Oma}[1]{\Omega_{#1} }
\nc{\zp}{{z'} }
\nc{\dps}{[ d \psi ] }

\begin{document}

\title{Localization
 for Nonabelian Group Actions}

\author{Lisa C. Jeffrey \\
Downing College \\
Cambridge CB2 1DQ, UK
\thanks{Address after 1 September 1993: Mathematics Department,
Princeton University, Princeton, NJ 08540, USA }  \\and\\
Frances C. Kirwan \\ Balliol College \\ Oxford OX1 3BJ, UK}

\date{February 1994}

\maketitle
\begin{abstract}

Suppose $X$ is a compact symplectic manifold
acted on by a compact Lie group $K$ (which may be
nonabelian) in a Hamiltonian
fashion, with moment map $\mu: X \to {\rm Lie}(K)^*$
and Marsden-Weinstein reduction $\xred
= \mu^{-1}(0)/K$.
There is then a
natural  surjective map $\kappa_0$ from the equivariant
cohomology $H^*_K(X) $ of $X$ to the cohomology
$H^*(\xred)$. In this paper we prove a formula
(Theorem 8.1, the  residue formula)
 for the evaluation on the fundamental class of
$\xred$ of any $\eta_0 \in H^*(\xred)$ whose
degree is the dimension of $\xred$, provided that $0$ is
a regular value of the moment map $\mu$ on $X$.
This formula is given in terms of any class $\eta \in H^*_K(X)$
for which $\kappa_0(\eta ) = \eta_0$, and involves
the restriction of $\eta$ to
$K$-orbits $KF$ of components $F \subset X$ of the fixed point
set  of a chosen maximal torus $T \subset K$.
Since $\kappa_0$ is
surjective,
in principle the residue formula enables one to determine generators
and relations for the cohomology ring $H^*(\xred)$, in
terms of generators and relations for $\hk(X)$.
 There are two main ingredients in the proof of our formula:
one is the
localization theorem  \cite{abmm,BV1}
for equivariant cohomology of manifolds acted on
by compact abelian groups, while the other is  the
equivariant normal form for
the symplectic form near the zero locus of the moment map.

We also make use of  the techniques appearing
in our proof of the residue formula
to give a new proof of the nonabelian localization formula of
Witten (\cite{tdg}, Section 2) for Hamiltonian actions of compact groups $K$
on symplectic manifolds $X$; this theorem expresses $\eta_0[\xred]$
in terms of certain integrals over $X$.

\end{abstract}
\renorm
\section{Introduction}

Suppose $X$ is a compact oriented manifold acted on by a compact
connected Lie group $K$ of dimension $s$; one may
then define the equivariant
cohomology $\hk(X)$.
Throughout this paper we shall consider only cohomology
with complex coefficients.
If $X$ is a symplectic manifold
with symplectic form $\om$ and the action of $K$ is Hamiltonian
(in other words, there is a moment map $\mu: X \to \lieks$), then we
may form the symplectic quotient $\xred = \mu^{-1}(0)/K$. The restriction
map $i_0: X \to \mu^{-1}(0)$ gives a ring homomorphism
$i_0^*: \hk(X) \to \hk (\mu^{-1}(0) ) $. Using Morse theory and the gradient
flow of the function $|\mu|^2: X \to \RR$, it is proved
in  \cite{Ki1}
that the map
$i_0^*$ is surjective.

Suppose in addition that $0$ is a {\em regular value} of the
moment map $\mu$. This assumption is equivalent to the assumption
that the stabilizer $K_x$ of $x$ under
 the action of $K$ on
$X$ is {\em finite } for every $x \in \mu^{-1} (0)$, and
it implies that $\xred$ is an orbifold, or $V$-manifold, which inherits
a symplectic form $\omega_0$ from the symplectic form
$\omega$ on $X$. In this situation there is a canonical isomorphism
$\pi_0^*: H^*(\mu^{-1}(0)/K) \to \hk (\zloc)$.\footnote{This
isomorphism is induced by the map $\pi_0: \zloc \times_K EK \to \zloc/K$.
Recall that we are only considering cohomology with complex coefficients.}
Hence we have a surjective ring homomorphism
\beq \label{0.1}
\k_0 = (\pi_0^*)^{-1} \circ i_0^*: \hk(X) \to
H^*(\xred). \eeq

Henceforth, if $\eta \in \hk(X)$ we shall denote $\k_0(\eta)$ by
$\eta_0$.
Previous work \cite{Brion,Ki2} on determining the ring structure
of $\hk(X)$ has presented methods which in some situations
permit the direct determination of the kernel of the
map $\k_0$, and hence  of generators and relations
in $H^*(\xred)$ in terms of generators
and relations in $\hk(X)$.  (Note that  the generators of $\hk(X)$
give generators of $H^*(\xred)$ via the surjective  map $\k_0$, and also that
generators of $H^*(BK)$ together with extensions to
$\hk(X)$ of generators of $H^*(X)$ give generators of
$\hk(X)$ because the spectral sequence
of the fibration
$X \times_K EK \to BK$ degenerates \cite{Ki1}.)
 Here we present
an alternative approach to determining the ring structure of
$H^*(\xred)$ when $0$ is a regular value of $\mu$, which complements
the results obtained by directly studying the kernel
of $\k_0$. Our approach is based on the observation that since
$H^*(\xred) $ satisfies Poincar\'e duality, a class
$\eta \in \hk(X)$ is in the kernel of $\k_0$ if and only if for all
$\z \in \hk(X)$ we have
\beq \eta_0  \zeta_0 [\xred]
= (\eta \zeta)_0 [\xred] = 0 . \eeq
Hence to determine the kernel of $\k_0$ (in other words the
relations in the ring $H^*(\xred)$) given
the ring structure of $\hk(X)$, it suffices to know
the {\em intersection pairings}, in other words the evaluations
on the fundamental class $[\xred]$ of all possible classes $
\xi_0 = \k_0 (\xi)$.
In principle the intersection pairings thus determine generators
and relations for the cohomology ring $H^*(\xred)$, given
generators and relations for $\hk(X)$.

 There is a natural
pushforward map $\pisk: \hk(X) \to \hk $ $ = \hk({\rm pt} )$
$ \cong S(\lieks)^K$, where we have identified
$\hk$ with the space of $K$-invariant polynomials on the
Lie algebra $\liek$. This map can be thought of as integration
over $X$ and will sometimes be denoted by
$\int_X$.
If $T$ is a compact {\em abelian} group
(i.e. a torus) and $\zeta \in \hht(X)$,
there is a formula\footnote{Atiyah and Bott
\cite{abmm} give
a cohomological proof of this formula, which was first proved by
Berline and Vergne \cite{BV1}.} (the {\em abelian localization
theorem})
for $\pist \zeta$ in terms of
the restriction of $\zeta$ to the components of the
fixed point set for the action of $T$. In particular, for a general
compact Lie group $K$ with maximal torus $T$ there
is a canonical map
$\tau_X: \hk(X) \to \hht(X)$, and we may apply the
abelian localization theorem to $\tau_X (\zeta)$ where
$\zeta \in \hk(X)$.

In terms of the components $F$ of the fixed point set
of $T$ on $X$, we obtain a
formula (the residue formula, Theorem \ref{t8.1})
for the evaluation
 of a class  $\eta_0 \in H^*(\xred) $ on the
fundamental class $[\xred]$, when $\eta_0$ comes
from a class $\eta \in \hk(X)$. There are two main
ingredients in the proof
of Theorem \ref{t8.1}. One is the abelian localization
theorem \cite{abmm,BV1},  while
the other is an equivariant  normal
form for $\om$  in a neighbourhood of $\zloc$,
given in \cite{STP}
as a consequence of the coisotropic embedding theorem.
The result is the following:

\noindent{\bf Theorem 8.1} {\em
Let $\eta \in \hk(X)$ induce $\eta_0 \in H^*(\xred)$.
Then we have
$$
\eta_0 e^{i\omega_0}  [\xred]  =
 \frac{(-1)^{n_+} }{(2 \pi)^{s-l}  |W| \vol(T) }
\treso
\Biggl (     \: \nusym^2 (\psi)
 \sum_{F \in \calf} e^{i \mu_T(F) (\psi) }
\int_F \frac{i_F^* (\eta(\psi) e^{i \omega} )  }{e_F(\mar  \psi) }
[d \psi] \Biggr ). $$
In this formula,  $n_+$ is the number of
positive roots of $K$,  and
$\nusym(\psi) = \prod_{\g > 0} \g(\psi)$ is the product of the
 positive roots, while
$\calf$  is the set of components of the fixed point
set of the maximal torus $T$ on $X$. If $F \in
\calf$ then $i_F$ is the inclusion of $F$ in $ X$ and
$e_F$ is the equivariant Euler class of the normal
bundle to $F$ in $X$.}

\noindent Here, via the {\em Cartan model},
the class $\tau_X ( \eta )\in \hht(X)$
has been identified with a family of differential
forms $\eta (\psi)$ on $X$
parametrized by $\psi \in \liet$. The definition   of the
residue map $\treso$ (whose domain is a suitable class of
meromorphic differential forms on
$\liet \otimes \CC$)
 will be given in Section
8 (Definition \ref{d8.5n}). It is a linear map, but
in order to apply it to the individual terms
in the statement of Theorem 8.1 some
choices must be made. The choices do not affect the
residue of the whole sum. When $\liet$ has dimension one the
formula becomes

$$
\eta_0 e^{i\omega_0}  [\xred]  =
 -\frac{1}{2}
{\rm Res}_0
\Biggl ( \psi^2      \:
 \sum_{F \in \calf_+} e^{i \mu_T(F) (\psi) }
\int_F \frac{i_F^* (\eta(\psi) e^{i \omega} )  }{e_F(\mar  \psi) }
 \Biggr ), $$
where ${\rm Res}_0$ denotes the coefficient of
$1/\psi$, and $\calf_+$ is the subset of the
fixed point set of $T = U(1)$ consisting of those
components $F$ of the $T$ fixed point set for which
$\mu_T(F) > 0 $.

We note that if $\dim \eta_0 = \dim \xred$ then the
left hand side of the equation in Theorem \ref{t8.1}
is just $\eta_0 [\xred]$. More generally one may obtain
a formula for $\eta_0[\xred]$ by replacing the
symplectic form $\omega$ by $\delta \omega$
(where $\delta > 0 $ is a small parameter), and taking
the limit as $\delta \to 0$. This  has the
effect of replacing the moment map $\mu$ by $\delta \mu$.
In the limit $\delta \to 0$, the Residue Formula (Theorem \ref{t8.1})
becomes a sum of terms corresponding to the components $F$ of the
fixed point set, where the term corresponding to $F$ is (up to a constant)
the residue (in the sense of Section 8) of
$\nusym^2(\psi) \int_F {i_F^* (\eta(\psi) )  }/{e_F(\mar  \psi) } $, and
the only role played by the symplectic form and the moment
map is in determining which
 $F$
give a nonzero
contribution to the residue of the sum and the signs with which
individual terms enter.

Results for the case when $K = S^1$, which are
related to our Theorem \ref{t8.1}, may
be found in the papers of Kalkman \cite{Kalk} and
Wu \cite{wu}.

Witten in Section 2 of \cite{tdg} gives a related result, the
{\em nonabelian localization theorem}, which also interprets
evaluations $\eta_0[\xred]$ of classes on the
fundamental class $[\xred]$ in terms of appropriate
data on $X$.  For $\e > 0$ and $\z \in \hk(X)$,
he defines\footnote{The normalization of the measure in
$\ie(\z)$ will be described at the beginning of
Section 3. As above, $\pis: \hk(X) \to \hk \cong S(\lieks)^K$ is the
natural pushforward map, where $S(\lieks)^K$ is the space of
$K$-invariant polynomials on $\liek$. }
\beq \label{1.1}
\ie(\zeta) = \frac{1}{(2 \pi i)^s \vol K}
\intk [d \phi] e^{- \e \inpr{\phi, \phi}/2 }
\pisk \z (\phi) \eeq
(where $\inpr{\cdot,\cdot} $ is a fixed invariant inner
product on $\liek$, which we shall use throughout to identify
$\lieks$ with $\liek$)
 and expresses it as a sum
of local contributions.

Witten's theorem tells us  that just as $\pisk \z$ would have
contributions  from
the components of the fixed point set of $K$ if $K $ were
abelian, the quantity $\ie(\z)$ (if $K$ is not necessarily
abelian) reduces to a sum of integrals
localized around the critical set of the function
$\rho = |\mu|^2,$
i.e. the set of points $x$ where $(d |\mu|^2)_x = 0$.
(Of course $d |\mu|^2 = 2\inpr{\mu, d \mu}$, so the
fixed point set of the $K$ action, where $d \mu = 0$,
is a subset of the critical set of $d|\mu|^2$.)
More precisely the critical set of $\rho = |\mu|^2$
can be expressed as a disjoint union of closed subsets
$C_\b$ of $X$ indexed by a finite subset
$\calb$ of the Lie algebra $\liet$ of the maximal
torus $T$ of $K$ which is explicitly
known in terms of the moment map
$\mu_T$ for the action of $T$ on $X$ \cite{Ki1}. If
$\beta \in \calb$ then the critical subset
$C_\beta$ is of the form $C_\beta = K (Z_\beta \cap
\mu^{-1}(\beta))$ where $Z_\beta$ is a union of
connected components of the fixed point set of the
subtorus of $T$ generated by $\beta$. The subset
$\mu^{-1}(0)$ on which $\rho = |\mu|^2$ takes its
minimum value is $C_0$.
There is a natural map\footnote{This map
is induced by the projection ${\rm pr}:\zloc \to {\rm pt}.$}
${\rm pr}^*: \hk \to \hk(\zloc)$
  so that the distinguished  class $f(\phi) = -\inpr{\phi, \phi}/2$
in $H^4_K$
gives rise to a distinguished class
$\T \in H^4(\zloc/K) \cong H^4_K(\zloc).$
Witten's result can then be expressed in the form
\begin{theorem} \label{t1.1}
$$ \ie(\z) = \z_0e^{\e \T}  [\xred]  + \sum_{\b \in
\calb - \{0\} }  \int_{U_\b}
{\z'}_\b. $$
Here,
the $U_\b$ are
open neighbourhoods in $X$ of
 the nonminimal critical subsets $C_\b$ of the function
$\rho$.
The $\zeb$ are
certain differential forms on $U_\b$ obtained from $\z$.
\end{theorem}

In the special case $\z = \eta \exp \iins \bar{\om} $ (where $\bar{\om}(\phi)
 = \om + \evab{\phi}{ \mu} $ is the standard extension of the symplectic
form $\om$
to an element of $H^2_K(X)$, and
$\eta$ has polynomial dependence\footnote{The equivariant
cohomology $\hk(X) $ is defined to consist of classes which have
polynomial dependence on the generators of $\hk$, but we shall
also make use of formal classes such as $\exp  \iins \bom$ which
are formal power series in these generators.}
on the generators of $\hk$),
 Witten's results give us   the following
 estimate
on the growth of the terms $\int_{U_\b} \zeb$
as $\e \to 0$:
\begin{theorem} \label{t1.2}
Suppose $\z = \eta \exp i \bar{\om}$ for some $\eta \in \hk(X)$. If
$\beta \in \calb - \{0\}$ then
$\int_{U_\b} \zeb = e^{- \rho_\b/{2 \e} } \:
h_\b(\e)$, where $\rho_\b = |\beta|^2$
 is the value of $|\mu|^2$
on the critical set $C_\b$ and $|h_\b(\e)|$ is bounded by a polynomial
in $\e^{-1}$.
\end{theorem}

Thus one should think of $\e > 0$ as a small parameter,
and one may use the asymptotics of the integral $\ie$
over $X$ to calculate the intersection pairings
$\eta_0 e^{\e \T} e^{i\om_0} [\xred],$ since the  terms in
Theorem \ref{t1.2} corresponding to the other critical subsets
of $\rho$ vanish exponentially fast as $\e \to 0$.
Notice that when $\z = \exp i\bom$, the vanishing of
$\mu$ on $\zloc$ means that $\z_0 = \exp i\om_0$, where
$\omega_0$ is the symplectic form induced
by $\omega$ on $\xred = \mu^{-1}(0)/K$.

In this paper we shall give a proof of a variant of
Theorems \ref{t1.1} and \ref{t1.2}, for the case
$\z = \eta \exp i \bom  $ where $\eta \in \hk(X)$.
 Before outlining our proof, it will
be useful to briefly recall Witten's argument.
Witten introduces a $K$-invariant 1-form $\l$ on
$X$, and shows that $\ie(\z) = \ie(\z \exp s D \l)$,
where $D $ is the differential in equivariant
cohomology and $s \in \RR^+$. He then does the
integral over $\phi \in \liek$ and shows that in the limit
as $s \to \infty$,
this
integral vanishes  over any region of
$X$  where $\l(\va)  \ne 0$ for at least one of the   vector fields
$\va, j = 1, \dots, s$ given by the infinitesimal action of a basis of
$\liek$ on $X$ indexed by $j$.
 Thus, after integrating
over $\phi \in \liek$, the limit as $s \to \infty$
of $\ie(\z)$ reduces to a sum of contributions
from sets where $\l(\va) = 0$ for all the $\va$.

In our case, when $X$ is a symplectic manifold and the action
of $K$ is Hamiltonian, Witten chooses $\l(Y) = d|\mu|^2(JY)$,
where $J$ is a $K$-invariant almost complex structure on
$X$. Thus $\l(\va)(x) = 0$ for all $j$
if and only if $(d|\mu|^2)_x = 0$,
so  $\ie(\z)$ reduces to a sum of contributions from the critical sets
of $\rho = |\mu|^2$.
Further, he obtains the contribution from $ \zloc$
as $e^{\e \T} \z_0 [\xred] $. If $\z = \eta e^{\iins \bom}  $
he also obtains the estimates in Theorem \ref{t1.2}
on the contributions from the neighbourhoods $U_\b$.

In general,
the contributions to the localization theorem depend
 on the choice of $\l$. In the symplectic
case, with $\l = J d |\mu|^2$, the contribution
from $\mu^{-1}(0)$ is canonical but the contributions
from the other critical sets $C_\b$ depend in principle
on the choice of $J$. Further, the properties of these
other terms are difficult to study. Ideally they
should reduce to integrals over the critical sets
$C_\b$, and indeed when proving  Theorem \ref{t1.2} Witten
makes   the assumption (before (2.52))  that the $C_\b$
are nondegenerate critical manifolds in the sense of Bott
\cite{bnd}.
In general the $C_\b$ are not manifolds; and
even when  they are manifolds, they are not necessarily nondegenerate. They
 satisfy only a weaker condition called  {\em minimal
degeneracy} \cite{Ki1}.\footnote{However, minimal degeneracy
may be sufficient for Witten's argument.} We shall treat the
integrals over neighbourhoods of the $C_\b$ in a future paper.

In the case when $X$ is a symplectic manifold and
$\z = \eta \exp \bom$ for any $\eta \in
\hk(X)$,
we have been able to use our methods
to  prove a variant
of Theorems \ref{t1.1} and \ref{t1.2}
(see Theorems \ref{t4.1}, \ref{t4.3} and \ref{t7.1} below)
which
bypasses these analytical difficulties and reduces the result
to fairly well known results on Hamiltonian group actions on
symplectic manifolds. We assume that $0$ is
a regular value of $\mu$, or equivalently
that $K$ acts on $\mu^{-1}(0)$ with finite
stabilizers.\footnote{Witten assumes that $K$ acts
freely on $\mu^{-1}(0)$.}  By treating the pushforward
$\pisk \z$ as a function on $\liek$, we
may use the  abelian
localization  formula \cite{abmm,BV1} for the
pushforward  in equivariant cohomology of
torus actions to find an explicit expression for
$\pisk \z$ as a function on $\liek$.
Thus, analytical problems relating to
integrals  over neighbourhoods of $C_\b$ are circumvented,
and localization reduces to studying the image of the moment
map and the pushforward of the symplectic or Liouville
measure under the moment map.\footnote{If $K$ is
abelian, this pushforward measure  is equal  (at $\phi \in \liek$)
to Lebesgue measure multiplied by  a function which gives
the symplectic volume of the reduced space
$\mu^{-1}(\phi)/K$; this function is sometimes called the
 Duistermaat-Heckman polynomial \cite{DH}.}
Seen in this light, the nonabelian localization
theorem is a consequence of the same results that
underlie the residue formula:
the  abelian localization formula for torus
 actions \cite{abmm,BV1} and the normal
form for $\om$ in a neighbourhood of $\zloc$.

We now summarize the key steps in our proof. Having replaced
integrals over $X$ by integrals over $\liek$,
we observe that in turn these may be replaced by integrals over the
Lie algebra
$\liet$ of the maximal torus. Then, applying properties of
the Fourier transform, we rewrite $\ie$ as the
integral over
$\liets$ of a Gaussian $\gtsoe(y) $ $\sim e^{- |y|^2/(2 \e) } $
 multiplied by a
function $Q = D_\nusym R$ where $R$ is piecewise
polynomial and $D_\nusym$ is a differential operator on $\liets$:
\beq \label{0.2}\ie =i^{-s} \tintt \widetilde{\gsoe} (y) Q(y),  \eeq
where $s$ is the dimension of $K$.
 The function $Q$ is obtained by combining the
abelian localization theorem (Theorem \ref{t2.1})
with a result
(Proposition \ref{p3.5}) on Fourier transforms of a certain class of
functions which arise in the
 formula for the pushforward.

The
function $Q$ is smooth in a neighbourhood of the origin when
$0$ is a regular value of $\mu$:
thus there is a polynomial $Q_0
 = D_\nusym R_0 $ which is equal to $Q$ near
$0$. It turns out that the cohomological expression
$\eeth[\xred]$ is obtained as the integral over $\liets$
of a Gaussian multiplied not
by $Q$ but by the polynomial $Q_0$:
\beq \eeth [\xred] = i^{-s}  \tintt \widetilde{\gsoe} (y) Q_0(y). \eeq
This result follows from a normal form for $\om$ near
$\zloc$.\footnote{This
normal form is a key tool in the original proof \cite{DH} of
the Duistermaat-Heckman theorem; this theorem
motivated  the proof by Atiyah and Bott  \cite{abmm} of the
abelian localization theorem.}

To obtain our analogue of Witten's  estimate (Theorem
\ref{t1.2}) for the asymptotics of
$\ie - \eeth [\xred]$ as $\e \to 0$,
we then write
\beq  \ie - \eeth [\xred] =
i^{-s} \tintt \widetilde{\gsoe}(y)
D_\nusym (R - R_0) (y). \eeq
Here, $R - R_0$ is piecewise polynomial and
supported {\em away} from $0$.
By studying
the minimum distances from $0$
in the support of $R - R_0$
we obtain an estimate (Theorem \ref{t4.1}) similar to
Witten's estimate (Theorem \ref{t1.2}). In our estimate,
the terms in the sum are indexed by
the set $\calb - \{0\}$; however, our estimate is
weaker than Witten's estimate since some of the subsets
$C_\beta$ indexed by
$\beta \in \calb - \{0\}$ (which a priori
contribute to our sum\footnote{In a future  paper
we hope to prove that the nonzero contributions
to our estimate (Theorem \ref{t4.1}) come only
from those $|\beta|^2$ which are nonzero critical values
of $|\mu|^2$.}
)
may be empty in which case $\rho_\beta
= |\beta|^2$ may not be a critical value
of $|\mu|^2$.

To summarize, the following related quantities appear in this paper:
\begin{enumerate}
\item The cohomological quantity $\eta_0 e^{i \om_0}
[\xred]$.
\item The integral $\cali^\epsilon $ (\ref{1.1}) coming from
the pushforward of an equivariant cohomology class
$\eta \in \hk(X)$  to $\hk$.
\item Sums of terms of the form
$$\int_F \frac{i_F^* \eta}{e_F}$$
where
$F$ is a connected component of
the  fixed point set  of the maximal torus $T$ acting
on $X$.
Such sums appear after mapping $\eta \in
\hk(X)$ into $\hht(X) $ and then applying the abelian
localization theorem.
\end{enumerate}
Witten's work relates (1) and (2), while our
Theorem \ref{t8.1} relates (1) and (3).

This paper is organized as follows. Section 2 contains
background material on equivariant cohomology
and the  abelian localization formula.
In Section 3 we collect a number of preliminary results
which we use in Section 4 to reduce our integral
$\ie$ to an integral over $\liets$ of a piecewise
polynomial function multiplied by a Gaussian.
Section 4 also contains the statement of two of
our main results,
Theorems \ref{t4.1} and \ref{t4.3}; Theorem
\ref{t4.3} is proved in Section 5, and Theorem \ref{t4.1}
in Section 6. In Section 7, Theorems \ref{t4.1}
and \ref{t4.3} (which are for the case $\z = \exp \bom$)
are extended to the case $\z = \eta \exp \bom $ for
$\eta \in \hk(X)$: the result is
Theorem \ref{t7.1}.  Finally, in Section 8 we
prove the  residue
 formula (Theorem \ref{t8.1}) for  the evaluation of cohomology
classes from $\hk(X)$ on the fundamental class of $\xred$, and
in Section 9 we
apply it when $K = SU(2)$ to specific examples.
This formula may be related to an unpublished
formula due to Donaldson.

In  future papers we shall treat the case when
$\xred$ is singular using intersection homology; we shall also
apply the nonabelian localization formula to moduli spaces of
bundles over Riemann surfaces regarded as finite
dimensional symplectic quotients, in singular as
well as nonsingular cases.

\noindent{\em Acknowledgement:}
We are most grateful to H. Duistermaat and M. Vergne for
helpful suggestions and  careful readings of the paper.

In addition, one
 of us (L.C.J.) wishes
to thank E. Lerman and E. Prato for explaining their work,
and also to thank E. Witten for discussions about
the nonabelian localization formula while the paper
\cite{tdg} was being written.

\renorm
\section{Equivariant cohomology and pushforwards}

In this section we recall the localization formula for torus
actions (Theorem
\ref{t2.1}) and express it in a form convenient
for our later use (Lemma \ref{l2.2}).

Let $X$ be a compact manifold equipped with the action of
a compact Lie group $K$ of dimension $s$ with
maximal torus $T$ of dimension $l$. We denote the Lie algebras
of $K$ and $T$ by $\liek$ and $\liet$ respectively, and
the Weyl group by $W$.
We assume an  invariant inner product
$\inpr{\cdot, \cdot} $ on  $\liek$ has
been chosen (for example, the Killing form): we shall use
this to identify $\liek$ with its dual. The orthocomplement
of $\liet$ in $\liek$ will be denoted $\lietp$.

Throughout this paper all cohomology groups are assumed to have
coefficients in
the field $\CC$.
The $K$-equivariant cohomology of a point is $\hk = H^*(BK)$,
and similarly the $T$-equivariant cohomology is $\hht
= H^*(BT)$. We  identify $\hk $  with
$S(\liek^*)^K$, the $K$-invariant polynomial functions on $\liek$,
and $\hht$ with $S(\liets)$.
Hence we have a bijective map  (obtained from the restriction from
$\lieks$ to $\liets$) which identifies $\hk$ with the subset of
$\hht$ fixed by the action of the Weyl group $W$:
\beq \label{2.1}
\hk \cong S(\lieks)^K
\cong S(\liets)^W \subseteq S(\liets) \cong \hht \eeq
This  natural map $\hk \to \hht$ will be denoted $\tau$ or
$\tau_X$.
We shall use the symbol $\phi$ to denote a point in
$\liek$, and $\psi$ to denote a
point in $\liet$. For $f \in \hk$ we shall write
$f = f (\phi)$ as a function of $\phi$.

The $K$-equivariant cohomology of $X$ is the cohomology
of a certain chain complex (see, for instance, Chapter 7 of
\cite{BGV} or Section 5
of \cite{MQ}; the construction is due to Cartan
\cite{cartan}) which
can be expressed  as
\beq \label{2.0}
\Om^*_K(X) = \Bigl ( S(\lieks) \otimes \Om^*(X)   \Bigr )^K \eeq
(where $\Om^*(X)$ denotes   differential forms on $X$).
An element in $\Om^*_K(X)$ may be thought of as a
$K$-equivariant polynomial function
from $\liek$ to $\Om^*(X)$. 
For
 $ \a \in \Om^*(X)$ and $f \in S(\lieks)$,
we write $(\a \otimes f) (\phi) = f(\phi) \a$.
In this notation, the differential $D$ on the
complex $\Om^*_K(X)$ is then defined  by\footnote{This definition and
the definition (\ref{2.0''}) of the extension $\bom(\phi)
= \omega +  \mu(\phi)$
of the
symplectic form $\omega$ to an equivariant cohomology
class are different from the conventions
used by Witten \cite{tdg}: a factor $\phi$ appears in the
our definitions where $i \phi$ appears in Witten's definition.
In other words Witten's definition is
$D(\a \otimes f) (\phi) =
f(\phi) (d \a -  i \iota_{ \tilde{\phi} } \a)$
and $\bom(\phi)
= \omega +  i \mu(\phi)$.
Witten makes this substitution
so that the oscillatory
integral $ \int_X \exp (\om + i \evab{\phi}{ \mu} ) $   will appear
as the integral of an equivariant cohomology class.}
\beq \label{2.0'}
D(\a \otimes f) (\phi) =
f(\phi) (d \a - \iota_{ \tilde{\phi} } \a)
=  \, f (\phi)  d \a
-  \sum_{j= 1}^s \phi_j f(\phi)  \: \iota_{\va} \a    . \eeq
Here, $\tilde{\phi}$ is the vector
field on $X$ given by the action
of $\phi \in \liek$, and  $\iota_{\tilde{\phi} }$ is the interior
product with the vector field $\tilde{\phi}$.
We have  introduced an orthonormal basis $\{ \hat{e}^j,
\: j = 1, \dots, s \} $ for $\liek$, and the $\phi_j$ $ \in \lieks$
are
simply the coordinate functions $\phi_j = \inpr{\hat{e}^j, \phi},
 $  while the $\va$ are the vector fields on $X$ generated by
the action of $\hat{e}^j$. The $\phi_j $ are assigned
degree $2$, so that the differential $D$ increases degrees by $1$.

One may define the pushforward
$\pis^K: \hk(X) \to \hk$, which
corresponds to integration over the fibre
of the map $X \times_K EK \to BK $ (see Section 2
of \cite{abmm}).  The pushforward  satisfies
$\pis^T = \tau \circ \pisk$. Because of this identification,
we shall usually simply write $\pis$ for $\pis^K$ or
$\pis^T$.
A localization formula for $\pist$ was given by
Berline and Vergne in \cite{BV1}; a more topological
proof of this formula is given in Section 3 of
\cite{abmm}.
\begin{theorem} \label{t2.1} {\bf \cite{BV1} }
If $\si \in \hht(X)$ and $\psi \in \liet$ then
$$(\pist \si) (\psi) = \sum_{F  \in \calf}
\int_F
\frac{i_F^* \si(\psi)  }{e_F( \psi) }. $$
Here we sum over the set $\calf$ of components $F$ of the fixed point
set of $T$, and $e_F $ is the $T$-equivariant
Euler class of the normal bundle of $F$; this Euler class
is an element of  $H^*_T(F) \cong H^*(F) \otimes \hht$,
as is $i_F^* \si$. The map $i_F: F \to X$ is the inclusion map.
The right hand side of the above expression is to be interpreted as
a rational function of $\psi$.
\end{theorem}

We shall now prove a lemma about the image of the pushforward,
which will be applied in Section 4.
\begin{lemma}\label{l2.2}
If $\si \in \hht(X) $ then $(\pist \si)(\psi)$ is a sum of terms
\beq \label{2.m0} (\pist \si)(\psi) =
\sum_{F \in \calf, \;\a \in \cala_F}\tau_{F,\a}  \eeq
such that each term $\tau_{F,\a}$  is of the form
\beq \label{2.m1}
\tau_{F,\a} =
\frac{\int_F c_{F, \a}(\psi)  }{ \efo(\mar \psi) \prod_j
\beta_{F,j}(\mar \psi)^{n_{F,j}(\a)} } \eeq
for some component $F$  of the fixed point set of the
$T$ action. Here, the $\beta_{F,j}$ are the weights of
the $T$ action on the normal bundle $\nu_F$, and
$\efo (\mar \psi) = \prod_j \beta_{F,j}(\mar \psi)$ is the product of all the
weights, while $n_{F,j}(\a)$ are some nonnegative integers. The class
$c_{F,\a}$ is in $H^*(F) \otimes H^*_T$, and is equal to
$i_F^* \si \in H^* (F)  \otimes H^*_T$ times some
characteristic
classes of subbundles of $\nu_F$.
\end{lemma}

\Proof
The normal bundle $\n_F$ to $F$ decomposes as a direct sum of
weight spaces  $\n_F  =
\oplus_{j = 1}^r \n_F^{(j)}$, on each of which
$T$ acts with weight $\bfj$. All these weights must be nonzero.
By passing to a split manifold if necessary (see section 21 of
\cite{BT}), we may assume without loss of generality
that the subbundle on which $T$ acts with a given weight decomposes
into a direct sum of $T$-invariant real subbundles
of rank $2$.
In other words, we may assume that the $
\nu_F^{(j)}$ are rank $2$ real bundles, and  the
$T$ action enables one to   identify them  in a
standard way with complex line bundles.

Then the equivariant
Euler class $e_F(\psi)$ is given for $\psi \in \liet$
 by
\beq \label{2.3}
e_F(\psi) = \prod_{j = 1}^r  \Bigl (c_1(\n_F^{(j)}) + \b_{F,j}(\psi)
\Bigr  ). \eeq
Thus we have
$$
\frac{1}{e_F (\mar \psi)}  = \frac{1}{\efo(\mar \psi) } \,
\prod_j (1 + \frac{c_1 (\nu_F^{(j)} )}{\beta_{F,j}(\mar \psi) } )^{-1} $$
\beq  \label{2.4}
= \frac{1}{\efo(\mar \psi) } \: \prod_j \sum_{r_j \ge 0 }
\: \, (-1)^{r_j} \Bigl ( \frac{c_1 (\nu_F^{(j)} )}
{\beta_{F,j}(\mar \psi) } \Bigr )^{r_j}
. \eeq
Here, $c_1(\nu_F^{(j)} ) \in H^2(F)$, so that
$c_1(\nu_F^{(j)} )/\b_{F,j}(\mar \psi)$ is {\em nilpotent} and the inverse
makes sense in $H^*(F) \otimes \CC(\psi_1, \dots, \psi_l)$,
where $\CC (\psi_1, \dots, \psi_l)$ denotes the complex
valued rational functions
on $\liet$.
$\square$

Let us now assume that $X$  is a symplectic manifold
and the action of $K$ is Hamiltonian with moment map
$\mu: X \to \lieks$. Denote by $\mu_T $
the moment map for the action of $T$ given by  the composition
of $\mu$ with the restriction map $\lieks \to  \liets$.
We shall be interested in one particular (formal)
equivariant cohomology
class $\si $, defined by
\beq \label{2.0''}\si(\phi) =
\exp  \iins \bom(\phi), \bla \bom(\phi)
 =  \omega +   \mu( \phi ) .  \eeq
 For this
class the localization formula gives
\beq \label{2.2}
 (\pist \si) (\psi) = \sum_F  \rf(\psi), \bla
\rf (\psi) =  \int_F \frac{  e^{i \mu_T(F)(\psi )
 } e^\om  } {e_F(\mar \psi) }. \eeq
(This formula does not require the fixed point set
of $T$ to consist of
 isolated fixed points.)

\noindent{\em Remark:} For any
$\eta \in \hk(X) $ the function $\pisk \eta \in \hk$ is
 a polynomial on $\liek$, and in particular is  smooth. However,
$\sigma = e^{\iins \bom} $ does not have
polynomial dependence on $\phi$.
Although it  is not immediately obvious
from the formula (\ref{2.2}), the
function $\pisk (\eta e^{\iins \bom}) $ is still  a
{ smooth}
function on $\liek$ (for any $\eta \in \hk(X)$ represented
by an element $\tilde{\eta} \in \Om^*_K (X)$): this
follows from its description as
$$ \pisk (\eta e^{\iins \bom}) (\phi) = \int_{x \in X} e^{i \om}
 \tilde{\eta} (\phi)
e^{i \mu(x) (\phi) }. $$

\renorm
\section{Preliminaries}

This section contains results which will be applied in the
next section to reduce the integral $\ie$ to an integral over $\liets$
of a Gaussian multiplied by a piecewise polynomial function. The
first, Lemma \ref{l3.1}, reduces integrals over $\liek$ to
integrals over $\liet$. Lemma \ref{l3.2} enables
us to replace the $L^2$ inner product of two functions by
the $L^2$ inner product of their Fourier transforms. Lemma
\ref{l3.3} relates Fourier transforms on $\liek$ to
Fourier transforms on $\liet$. Finally Proposition \ref{p3.5}
describes  certain functions
whose Fourier transforms are the terms
appearing
in the  localization formula (\ref{2.4}).

We would like to study a certain integral that arises out
of equivariant cohomology:
\beq \label{3.1}
\cali^\e = \frac{1}{(2 \pi i )^s\vol K} \int_{\phi \in \liek}
[d\phi] \, e^{- \e \inpr {\phi, \phi } /2 }  \int_X \si(\mar \phi). \eeq
Here, $\si \in \Om^*_K(X) $
(see (\ref{2.0})); we are mainly interested in the
class $\si$ defined by (\ref{2.0''}).
Also,
 $\e > 0 $ and we shall consider the behaviour of $\cali^\e$
as $\e \to 0^+$.
The measure $[d \phi]$ is a measure on $\liek$ which
corresponds to a choice of invariant metric on $\liek$
(for instance, the metric given by the Killing form):  such
a metric induces a  volume form on $K$, and
$\vol K$ is the integral of this volume  form over $K$. Thus
$[d \phi]/\vol K$ is independent of the choice of metric on
$\liek$. The metric also gives a  measure  $[d \psi] $
on $\liet$ and a volume form on $T$: it is
implicit in our notation that the measures  on $T$ and $\liet$
come from the same invariant metric as those on $K$ and $\liek$.

It will be convenient to recast integrals over $\liek$ in terms
of integrals over $\liet$. For this we use
  a function $\nusym: \liet \to \RR$,
satisfying $\nusym(w \psi) = (\det w) \nusym(\psi)$ for
all elements $w$ of the Weyl group  $ W$, and
defined by
\beq \label{3.2}
\nusym(\psi) = \prod_{\g > 0} \g(\psi), \eeq
where $\g$ runs over the positive roots.
Using the inner product to identify $\liet$ with
$\liets$, $\nusym$ also defines a function $\liets \to \RR$.
We  have
\begin{lemma} \label{l3.1}[Weyl Integration Formula]
If $f: \liek \to \RR$ is $K$-invariant, then
$$\int_{\phi \in \liek} f (\phi) [d \phi]
= \ck^{-1}  \intt f (\psi) \nusym(\psi)^2
[d \psi], $$
where $s$ and $l$ are the dimensions of $K$ and $T$, and
$ \ck =  { \wn \vol T}/\vol K $.
\end{lemma}
\Proof There is an orthonormal basis $\{X_\g, Y_\g | $
$\g $ a positive root$\}$ for $\lietp$ such that
$$[X_\g, \psi] =  \g(\psi) Y_\g, $$
$$ [Y_\g, \psi] = -  \g(\psi) X_\g$$
for all $\psi \in \liet$. The Riemannian volume form of the
coadjoint orbit
through $\psi \in \liet \cong \liets$ (with the metric
on the orbit pulled back from the metric on $\lieks$
induced by the inner product $\inpr{\cdot, \cdot}$)
evaluated on the tangent vectors $[X_\g, \psi]$ and $[ Y_\g, \psi]$
is thus $ \prod_{\g> 0} \g(\psi)^2$, while the volume
form
of the homogeneous space $K/T$ (induced by
the chosen metric on $\liek$) evaluated on the tangent vectors
corresponding to
$X_\g, Y_\g \in \liek $ is $1$. Hence
 the Riemannian volume of the  orbit
through $\psi \in \liet$ is $ \nusym(\psi)^2$ times the
volume of the homogeneous space $K/T$.
$\square$

\nc{\lamax}{\L^{\rm max} }
\nc{\dist}{\cald'}

It will be convenient also to work with the Fourier transform. Given
$f: \liek \to \RR$ we define $F_K f: \lieks \to \RR$,
$F_T f: \liets \to \RR$
by
\beq \label{3.3}
(F_K f) (z) = \frac{1}{(2 \pi)^{s/2} } \intk f (\phi)
e^{-i \evab{\phi}{z} }  \, [d \phi], \eeq
\beq \label{3.4}
(F_T f) (y) = \frac{1}{(2 \pi)^{l/2} } \intt f (\psi)
e^{-i \evab{\psi}{y} } \, [d \psi]. \eeq
More invariantly, the Fourier transform
is defined on a vector space $V$  of dimension $n$ with
dual space $V^*$ as a map
$F: \O^{\rm max}(V) \to \O^{\rm max} (V^*)$
where
$\O^{\rm max}  (V) = \L^{\rm max} (V^*) \otimes \dist (V)$,
$\L^{\rm max}(V^*) $
  is the top exterior power of
$V^*$ and $\dist (V) $ are the tempered distributions
on $V$ (see \cite{hor}).
Indeed for $z \in V^*$, $f \in \dist(V) $ and $u \in \lamax(V^*)$
we define
\beq \Bigl ( F (u \otimes f )\Bigr ) (z) = \frac{v}{(2 \pi)^{n/2}}
\int_{\phi \in V}  f(\phi) e^{- i z(\phi)}  u  , \eeq
where the element
$v\in \lamax (V) $ satisfies $u(v)  = 1 $ under the
natural pairing $\lamax(V^*) \cong \Bigl (\lamax (V) \Bigr)^*.$
(The normalization has been chosen so that
$$ \tilde{f} (\phi) = F (F \tilde{f}) (- \phi)  $$
for any $\tilde{f} \in \O^{\rm max}(V)$.)
For notational convenience
we shall often ignore this subtlety
and identify $\liek, \liet$ with $\lieks, \liets$
under the invariant inner product $\inpr{\cdot, \cdot}$:
we shall also suppress the exterior powers of
$\liek$ and $\liet$ and write $\fk: \dist(\liek) \to
\dist (\lieks)$, $\ft: \dist(\liet) \to \dist(\liets)$.
Further, although we shall work with functions whose definition
depends on the choice of the element  $[d \phi] \in \lamax(\lieks)$
(associated to the inner product), some of our   end results
do not depend on this choice\footnote{The statement of
Theorem \ref{t8.1},
for instance, does not depend on
the invariant
inner product $\inpr{\cdot, \cdot}$.}
and the use  of such
functions is just  a notational convenience.

A fundamental property of the Fourier transform
is that it preserves the $L^2$ inner product:
\begin{lemma} \label{l3.2} {\bf [Parseval's Theorem]}(\cite{hor},
Section 7.1)
If $f:  \liek \to  \CC$  is a tempered distribution
and $g: \liek \to \CC$ is a Schwartz function then
$F_K f: \liek \to \CC$ is also a tempered distribution and
$F_K g: \liek \to \CC$ a Schwartz function, and we have
$$ \intk \overline{g(\phi)}
 f (\phi) [d \phi] = \tintk \overline{(F_K g)(z) }
(F_K f)(z) \, [d z], $$
$$ \intt \overline{g(\psi)}
 f (\psi) [d \psi] = \tintt \overline{(F_Tg)(y) }
(F_T f)(y) \, [d y]. $$
\end{lemma}
We note also that if
$g_\e: \liek \to \RR$ is the
Gaussian defined by $g_\e(\phi) =  e^{-\e \inpr{\phi, \phi}/2  }
$, then
\beq \label{3.5} (F_K g_\e) (z) = \frac{1}{\e^{s/2} } e^{- \inpr{z,z}
/{2 \e}
} = \frac{1}{\e^{s/2} }  \gsoe(z),  \bla
(F_T g_\e  ) (y) = \frac{1}{\e^{l/2} } e^{- \inpr{y,y} /{2 \e}
} =  \frac{1}{\e^{l/2} } \gsoe(y). \eeq

We have also that
\begin{lemma} \label{l3.2'} The symplectic volume form
$d \Om^S_\phi$ at a point $\phi$ in the
orbit $K \cdot \psi$
  through $\psi \in \liek$
is related to the Riemannian volume form
$d \Om^R_\phi$ (induced by
the metric on $\liek$)
by
$$d \Om^R_\phi =  \nusym(\psi) d \Om^S_\phi.$$
\end{lemma}
\Proof The symplectic form is
  $K$-invariant
and is given at the point $\psi \in  \liet_+$  in the orbit
$K \cdot \psi$
(for $\xi, \eta \in \lietp$ giving rise to tangent
vectors $[\xi,\psi], $ $[\eta, \psi]$ to the orbit) by
\beq \label{3.7}
\om([\xi,\psi], [ \eta,\psi]) = \inpr{[ \xi,\psi], \eta }
= \inpr{\psi, [ \eta,\xi] }.  \eeq
In the notation of Lemma \ref{l3.1}, the symplectic
volume form evaluated on the tangent vectors
$[ X_\g,\psi], $ $[ Y_\g,\psi]$ is given by
$\prod_{\g > 0} \om ([ X_\g,\psi],  [ Y_\g,\psi] ).$
But $\om ([X_\g,\psi],  [ Y_\g,\psi] ) =  \inpr{[ X_\g,\psi], Y_\g} $
$ =  \g(\psi)$, from which (comparing with the proof of
Lemma \ref{l3.1}) the Lemma follows. $\square$

We shall use this Lemma to prove the following
Lemma  relating Fourier
transforms on $\liek$ to those on $\liet$:
\begin{lemma} \label{l3.3}
Let $f \in \dist(\liek)$ be $K$-invariant,
and let $\nusym$ be defined
by (\ref{3.2}). Then
$$\ft(\nusym f) = \nusym \fk(f)$$
as distributions on $\liet$.
\end{lemma}
\Proof $$(\fk f) (z) = \frac{1}{(2 \pi)^{s/2} } \intk e^{-i
\evab{\phi}{z} } \, f(\phi) [d \phi] $$
\beq \label{3.6}
= \frac{1} {(2 \pi)^{s/2} } \int_{\psi \in \liet_+}
\int [d \psi] f(\psi) \: \int_{\phi \in K \cdot \psi}
e^{- i \evab{\phi}{z} } d \Om^R_\phi \eeq
(where $\liet_+$ denotes the fundamental Weyl chamber).

We have from Lemma \ref{l3.2'} that
$d \Om^R_\phi =  \nusym(\psi) d \Om^S_\phi$.
Thus we have
\beq \label{3.8}(F_K f)(z) = (2 \pi)^{-s/2} \int_{\psi \in \liet_+}
[d \psi] f (\psi) \nusym(\psi) \int_{\phi \in K \cdot \psi}
e^{- i \evab{\phi}{ z} }  d \Om^S_\phi. \eeq

Now the integral over the coadjoint orbit may be computed
by the Duistermaat-Heckman theorem \cite{DH} applied to the
action of $T$ on the orbit $K \cdot \psi$
(or equivalently as a consequence
of the
abelian localization theorem, Theorem \ref{t2.1}): we have
(see e.g. \cite{BGV} Theorem 7.24)
\beq \label{3.9}
\int_{\phi \in K \cdot \psi} e^{- i \evab{\phi}{ z} }
\, d \Om^S_\phi =
\frac{(2 \pi)^{(s-l)/2} } {\nusym(z) } \sum_{w \in W} e^{- i \evab{w \psi}{ z}
}
(\det w). \eeq
This is a well known formula due originally to Harish-Chandra
(\cite{Harish}, Lemma 15).
(Notice that the  $z$ in (\ref{3.9}) plays the role of the
$\psi$ in Theorem \ref{t2.1}, while the $\psi$ in
(\ref{3.9}) specifies the orbit.)
Now since $\nusym(w \psi) = (\det w) \nusym(\psi), $ we may replace
the integral over $\liet_+$ by an integral over $\liet$: in other
words we have
\beq \label{3.10}
\nusym(z) (F_K f) (z) = (2 \pi)^{- l/2} \int_{\psi \in \liet}
[d \psi] \, f(\psi) \nusym(\psi) \eminevb{\psi}{ z} = F_T(\nusym f) (z).
\bla \square \eeq

\nc{\ims}{i^{-s} }

Applying this to the Gaussian
$\gse(\psi) =  e^{-\e \inpr{\psi, \psi}/2  }  $ we have
\begin{corollary} \label{c3.4}
$$
F_T(\gse \nusym) (y) = \frac{1}{\e^{s/2} } \nusym(y) e^{- \inpr{y,y} /{2 \e} }
= \frac{1}{\e^{s/2} } \nusym(y)  \gsoe(y).  $$
\end{corollary}

We shall also need
the following result which occurs
in the work of Guillemin, Lerman, Prato
and Sternberg
 concerning Fourier transforms of
a class of functions on
$\liets$. This result will be applied to functions
appearing in the abelian localization formula (\ref{2.4}).
\begin{prop} \label{p3.5}
{\bf (a)} (see \cite{JGP},  section 3.2 of \cite{GLS},
and \cite{GP})
Define $h(y) = H_{\bar{\beta} }(y + \tau)$ for some
$\tau \in \liet$,
where $H_{\bar{\beta} } (y) = \vol  \{ (s_1, \dots,
s_N): s_i \ge 0,  \onebl y = \sum_j s_j \beta_j \}$ for
some $N$-tuple   $\bar{\beta} = $
$\{ \b_1, \dots, \b_N \}$ , $\b_j \in \liets$,
such that
the $\b_j$ all lie in the interior of
 some half-space of $\liets$.
Thus $H_{\bar{\beta} } $ is a piecewise polynomial
function supported on the cone $C_\barb =
\{ \sum_j s_j \beta_j \, | \, s_j \ge 0 \}. $
Then the Fourier transform of $h$ is given for
$\psi$ in the complement of the hyperplanes
$\{ \psi \in \liet| {\beta_j}({ \psi}) = 0 \}$ by
the formula
\beq \label{3.10'} F_T h(\psi) = \frac{\epinev{\tau}{ \psi} }
{i^N \prod_{j = 1}^N \beta_j (\mar \psi) }.  \eeq

\noindent{\bf(b)} (see Section 2 of \cite{JGP}
and (2.15) of \cite{GP}) The function $H_{\bar{\b}} (y) $ is also given
as
$$H_{\bar{\b}} (y) = H_{\b_1} * H_{\b_2} * \dots
* H_{\b_r}, $$
where for $\b \in {\rm Hom}(\RR^l,\RR)$ we have  $H_\b = (i_\b)_* dt$,
i.e. $H_\b$ is the pushforward of the Euclidean
measure $dt$ on $\RR^+$ under the map $i_\b: \RR^+ \to \RR^l$
given by $i_\b(t) = \b t$.  Here, $*$ denotes convolution.

\noindent{\bf (c)} (see Section 2 of \cite{JGP}) The function
$H_\barb $ satisfies the differential equation
$$ \prod_{j = 1}^N \b_j (\partial/\partial y) H_\barb(y) =
\delta_0(y) $$
where $\delta_0$ is the Dirac delta distribution.

\noindent{\bf (d)} (see Proposition 2.6 of \cite{JGP})
The function $H_\barb$ is smooth at any
$y \in U_\barb$, where $U_\barb$ are the
points in $\liets$  which are not in any cone
spanned by a subset of $\{\b_1, \dots, \b_N \} $ containing
fewer than $l$ elements.

\end{prop}

\renorm
\section{Reduction to a piecewise polynomial function}

	In this section we apply the results stated in
Section 3 to reduce $\ie$ to the form given in Proposition \ref{p4.2},
as the integral over $y \in \liets$ of a
Gaussian $\gtsoe(y) = \gsoe(y)/((2 \pi)^s  \wn \vol( T) \e^{s/2} ) $
times a function
$Q(y)$ which is $D_\nusym R(y)$
where $D_\nusym$ is a differential operator
and $R$ is a piecewise polynomial function.
As a byproduct we obtain also a generalization of
Theorem 2.16 of \cite{GP}: we may relax the hypothesis of \cite{GP}
that
the action of the maximal torus $T$ have isolated fixed points.

 Two of our main results related to Witten's work in
\cite{tdg}
are Theorem \ref{t4.1} and \ref{t4.3}: these are stated in this section.
The  proof of Theorem \ref{t4.3} will be given in Section 5. It
tells us  that the cohomological contribution $\eeth[\xred]$
given in Witten's Theorem \ref{t1.2} for the zero locus
of the moment map is given by the integral
over $y \in \liets$ of $\gtsoe(y)$ times a polynomial
$Q_0(y)$ which is equal to $Q$ near $y = 0$.
Theorem \ref{t4.1} is our version of the asymptotic
estimates given in Witten's Theorem \ref{t1.1}, and will be
proved in Section 6 below. Theorem \ref{t7.1} in Section 7
extends Theorems \ref{t4.1} and \ref{t4.3} to more general
equivariant cohomology classes.

We assume throughout the rest of the
paper that $K$ acts on $X$ in a Hamiltonian fashion,
and that $0$ is a regular value of the moment map
$\mu$ for the $K$ action.
This is equivalent to the assumption that
 $K$ acts on $\mu^{-1}(0)$ with finite stabilizers, and it
 implies that $\mu^{-1}(0)$ is a smooth manifold.
 Under these
hypotheses, the space $\xred = \mu^{-1}(0)/K$ is a $V$-manifold or
orbifold
(see \cite{kaw}) and $P = \mu^{-1}(0) \to
\mu^{-1}(0)/K $ is a $V$-bundle: we  have a class
$\T \in H^4 (\xred) $ which represents the class
$-\inpr{\phi, \phi}/2 \in H^4_K(\mu^{-1}(0) ) $
$\cong H^4(\zloc/K)$,
and which is a four-dimensional characteristic class
of the bundle $P \to \xred$.

In \cite{Ki1} it is proved that the set of
critical points of the function $\rho = |\mu|^2: X \to \RR$
is a disjoint union of closed subsets $C_\b$ in
$X$ indexed by a finite subset $\calb$ of $\liet$.
In fact if $\liet_+$ is a fixed positive Weyl
chamber for $K$ in $\liet$ then $\b \in \calb$ if and only
if $\b \in \liet_+$ and $\b$ is the closest point to $0$
of the convex hull in $\liet$ of some nonempty
subset of the finite set $\{ \mu_T(F): F \in \calf\}$,
i.e. the image under $\mu_T$ of the set of fixed
points of $T$ in $X$. Moreover if $\b \in \calb$
then $$C_\b = K(Z_\b \cap \mu^{-1}(\b))$$
where $Z_\b$ is the union of those connected components
of the set of critical points of the function $\mu_\b$
defined by $\mu_\b (x) = \eva{\mu(x)}{  \b} $
on which $\mu_\b$ takes the value $|\b|^2$. Note that
$C_0 = \mu^{-1}(0) $ and in general the value
taken by the function $\rho = |\mu|^2$
on the critical set $C_\b$ is just $|\b|^2$.

We shall prove the following version
of Witten's nonabelian localization theorem, for
the integral $\ie$  defined in (\ref{3.1}) with
the
class $\si  = e^{\iins \bom} $ defined by (\ref{2.0''}):
\begin{theorem}
\label{t4.1}  For each $\b \in \calb$ let $\rho_\b = |\b|^2$
(this is   the  critical value
of the function $\rho = |\mu|^2: X \to \RR$ on the critical
set $C_\b$ when this set is nonempty).
Then
there exist functions $h_\beta: \RR^+ \to \RR $  such that
for some $N_\beta \ge 0$,
$\e^{N_\beta} h_\beta(\e)$ remains bounded as $\e \to 0^+$, and
for which
$$|\ie  - e^{\e \T} e^\om [\xred]  | \le
\sum_{\beta \in \calb - \{0\} }
 e^{- \rho_\beta  /{2 \e}  } \, h_\beta(\e). $$

\end{theorem}

\noindent{\em Remark:} The estimate given in Theorem
\ref{t4.1} is {\em weaker} than Witten's estimate
(Theorem 1.2) since $|\beta|^2 $ is not in fact a
critical value of $|\mu|^2 $ when $C_\beta$ is
empty.
Nevertheless in many interesting cases all the $C_\beta$
are nonempty, and our estimate then coincides with Witten's.

To prove this result, we shall rewrite $\ie$ using Lemma
\ref{l3.1}:
\beq \label{4.1} \ie = \frac{1}{(2 \pi i)^s  \wn \vol(T) }
 \: \intt [d \psi] \, \Bigl (
\gse(\psi) \: \nusym(\psi) \Bigr ) \:
\Bigl ( \: \nusym(\psi) (\pist \sigma ) (\psi) \Bigr ). \eeq
Now we apply Lemma \ref{l3.2} to get
\beq \label{4.2}
\ie = \frac{1}{(2 \pi i)^s \wn \vol(T) }
\tintt [dy] \Bigl ( F_T (g_\e \nusym) (y) \Bigr )
\Bigl ( F_T (\nusym \, \pist \si ) (y) \Bigr ). \eeq
Applying Corollary \ref{c3.4} we have
\beq \label{4.3}
\ie = \frac{1}{(2 \pi i)^s \wn \vol(T) \e^{s/2}  }
\tintt [dy] \nusym(y) e^{- \inpr{y,y}/{2 \e} } \,
F_T ( \nusym \pist \si) (y). \eeq

\nc{\bwd}{\lasub}
\nc{\cf}{C_F}
\nc{\cfl}{C_{F,\lasub } }
\nc{\dcf}{\check{C}_{F,\lasub} }

Following Guillemin, Lerman and Sternberg \cite{JGP},
we may  use the abelian localization formula (Theorem
\ref{t2.1})
to give a formula for $F_T ( \pist \si)$ where $\si = e^{\iins \bom}$,
and from it obtain a formula
for $F_T (\nusym \pist \si)$.
In terms of the notation of  Lemma \ref{l2.2}, we choose
a component $\bwd$ of the set $\cap_{F,j}  \Bigl \{
\psi \in \liet:$ $ \bfj(\psi) \ne 0 \Bigr \}, $
where $\b_{F, j}$ are the weights of the action
of $T$ on the normal bundle to a component $F$ of the
fixed point set.  Thus $\bwd$
is a cone in $\liet$. If we denote by $\cf = \cfl$
the component of $\cap_j \{ \psi \in \liet:
$ $ \bfj (\psi) \ne 0 \} $ containing $\bwd$, then
$\bwd = \cap_F \cfl$. Also, $\bfj \in \liets$ lies
in the dual cone $\dcf$ of $\cfl$: indeed, this dual
cone is simply the cone $\dcf = \{
\sum_j s_j \bfj: s_j \ge 0 \}. $
We then define
$\s_\fj = {\rm sign} \bfj(\xi) $ for any
$\xi \in \bwd$, and $\bfjw =
\s_\fj \bfj$. Then we set $k_F(\a) =
\sum_{j, \s_\fj = - 1} n_\fj(\a)$.

We define a function  $\hh: \liets \to \CC$
by
\beq \label{4.g}
\hh(y) = \sum_{F \in \calf} \sum_{\a  \in \cala_F}
(-1)^{k_F(\a) } H_{\bar{\g_F}(\a)  } (-y + \mu_T(F) )
\int_F (e^{i\om} \tilde{c}_{F,\a}) . \eeq
Here as before $\calf$ is the set of components
$F$ of the fixed point set of $T$ and
the $\cala_F$ are the indexing sets which appeared
in Lemma \ref{l2.2}. If $\a \in \cala_F$ then
$\bar{\g_F}(\a)$ consists of the elements
$\bfjw$ where each $\bfjw$ appears with multiplicity
$n_\fj(\a) $.
Then
$H_{\bar{\g_F}(\a) } $ is  as defined in Proposition \ref{p3.5}.
The $\tilde{c}_{F,\a} \in H^*(F) $ are related to
the $c_{F,\a} $  in (\ref{2.m1}), in that
\beq \label{4.004} c_{F, \a}(\psi)  = e^{i\om}  e^{i \eva{\mu_T(F)}{ \psi} }
\tilde{c}_{F,\a}.  \eeq

We then have the following Theorem, which in the case
when the action of $T$ has isolated fixed points is the main
theorem
of Section 3 of the  paper \cite{JGP}
of Guillemin, Lerman and Sternberg.
For the most part our proof is a direct extension
of the proof given in that paper; the major difference is in
the use of the abelian localization theorem rather than
stationary phase.

\begin{theorem} \label{t4.1'} The (piecewise polynomial)
function
$\hh$ given in (\ref{4.g}) is identical
to  the distribution $ \ggh  = F_T ( \pist e^{\iins \bom}) $.
\end{theorem}
\Proof We first
 apply the abelian localization formula (\ref{2.2})
to $\pist \si$.
 Then
the formula obtained from Lemma \ref{l2.2}  for
$\pist \si$ is
\beq \label{4.3p} \pist \si =
\sum_{F \in \calf} \rf, \bla \rf = \sum_{\a \in \cala_F} \tau_{F, \a} , \eeq
\beq \label{4.03p}\tau_{F, \a} =
 (-1)^{k_F(\a) }  \frac{ \epinevb{\psi}{ \mu_T(F)}
\int_F (e^{\iins\om} \tilde{c}_{F,\a} ) }
{\prod_j \Bigl ( \bfjw (\mar \psi) \Bigr )^{{n}_\fj(\a)} } , \eeq
where $\mu_T $, the moment map for the $T$ action, is simply
the projection of $\mu$ onto $\liet$
and the  $ \bfj  $ are the weights of the $T$ action. Recall that
the quantity $\bfjw$ is  $\bfj$   if $\bfj (\xi) > 0 $
and $- \bfj$ if $\bfj (\xi) < 0 $.
The class $\tilde{c}_{F, \a} \in H^*(F)$ is equal to
some characteristic classes of
subbundles of the normal bundle $\nu_F$; it is independent
of $\psi$.

Notice that each  $\bfjw$ $\in \liets$
lies in the half space
$\{ y  \in \liets \, | y( \xi) > 0 \}$, for any $\xi \in \bwd$.
The expression  (\ref{4.03p}) is hence
of the form  appearing on the right hand side
of (\ref{3.10'}), up to multiplication  by a
factor independent of $\psi$.
The conclusion of the proof goes, as in Section 3 of \cite{JGP},
by applying a lemma about distributions (see Appendix A of \cite{GP}):
\begin{lemma} \label{l4} Suppose $\gggh$ and $\hhh$ are
two tempered distributions on $\RR^l$ such that:
\begin{enumerate}
\item $F_T \gggh - F_T \hhh$  is supported on a finite union
of hyperplanes.
\item There is a half space $\{ y \, | \, \inpr{y,\z} > k_0 \}$
containing the support of $\hhh - \gggh$.
\end{enumerate}
Then $\gggh = \hhh$.
\end{lemma}
Here we apply the lemma to $\hhh$ as given in
(\ref{4.g}) and $\gggh = F_T ( \pist e^{\iins \bom} ) $.
The first hypothesis  is satisfied  because we know from Proposition
\ref{p3.5} that $F_T \hh$  is given by the formula (\ref{4.3p})
on the complement of the hyperplanes
$\{ \psi \, | \, \bfjw ( \psi)  = 0 \}$; but this
is just the formula for $F_T \ggh = \pist e^{\iins \bom}$.
Further, $\hh$  is supported in a half space since all the
weights $\bfjw$ satisfy $\bfjw(\xi) > 0$ for any $\xi \in \bwd$, while
the support of $\ggh$ is
contained in the compact set $\mu_T(X)$  (see Section 5 below).
Therefore the support of $H-G$ is contained
in a half space of the form $\{ y: \inpr{y,\xi} > k_0 \}$
for some $k_0$.
This completes the proof of Theorem \ref{t4.1'}. $\square$

Define
\beq \label{4.3pp}
R(y) =  F_T ( \pist \si) (y), \eeq
where $\si = \exp \iins \bom$.
Then
 (\ref{4.3}) and (\ref{4.g}) give us  the following
\begin{prop} \label{p4.2}
The function  $R$ is a piecewise polynomial function supported on
cones each of which has apex at $\mu_T(F)$ for some component
$F$ of the fixed point set of $T$.
Let $Q$ be the distribution defined by
$$ Q(y) = \nusym(y) D_\nusym R(y) $$
where
the differential operator $D_\nusym$ is given by
$$D_\nusym = \prod_{\g> 0} (i \g(\partial/\partial y) )$$
and $\g$ runs over the positive roots (cf. (\ref{3.2})).
Then
\beq \label{4.3ppp}
\ie = \frac{1}{(2 \pi i )^s \wn \vol( T) \e^{s/2} }
\tintt [dy] e^{- \inpr{y,y}/{2 \e}  }  Q(y). \eeq
\end{prop}
\Proof It only remains
to note that  $\ft (\nusym \pist \si) = D_\nusym \ft (\pist \si), $
where $D_\nusym$ is defined above. $\square$

\noindent{\em Remark:} The formulas for $Q(y)$ obtained
from (\ref{4.3p}) will in general be different for
different choices of $\bwd$.

In addition, certain formulas simplify if we impose
the additional assumption that at any point $x$ in a component
$F$ of the fixed point set of the $T$ action,
the orthocomplement $\lietp$ of $\liet$ in $\liek$
injects into the tangent space $T_x X$ under the infinitesimal
action of $K$: in other words, that the stabilizer ${\rm Stab} (x)$
of $x$ is such that ${\rm Stab}(x)/T$ is a finite group.
Under this  additional hypothesis,
we may indeed prove a somewhat stronger result.
 Notice that by Lemma \ref{l2.2},
each term (corresponding to a component $F$ in the fixed point
set of $T$) in the localization
formula for $\pist \si$ has a factor $\efo(\mar \psi)$ in the denominator.
Now
for $x \in F$, the fibre $(\nu_F)_x$ over $x$  of the normal bundle
to $F$
will contain $ \lietp \cdot x$, the image of $\lietp$
under the infinitesimal action of $K$. Under the additional assumption
that $\lietp$ injects into $T_x X$,
the set of weights for $\nu_F$ contains for each root $\g > 0$
either the root $\g$ or the root $-  \g$: in other words,
$\efo (\mar \psi)$ is divisible by $\nusym(\psi)$. Thus from Lemma
\ref{l2.2} we obtain
the formula for
$\nusym(\psi) \pist \si$:
\beq \label{4.3pinj}
\nusym(\psi) \pist \si = \sum_{F \in \calf, \a \in \cala_F}
\ttfa, \eeq
\beq \label{4.03pinj} \ttfa =
 (-1)^{k_F(\a) }  \frac{ \epinevb{\psi}{ \mu_T(F)}
 }
{\prod_j \Bigl ( \bfjw (\mar \psi) \Bigr )^{\tilde{n}_\fj(\a)} }
\int_F ( e^{i \om} \tilde{c}_{F,\a} ) . \eeq
Here,  the notation is as in (\ref{4.3p})
and (\ref{4.03p}) except that
$\tilde{n}_\fj(\a) = n_\fj(\a)  $ if $\bfj$ is not a
root, while   $\tilde{n}_\fj(\a) = n_\fj(\a) - 1$ if $\bfj$ is a root.

 We may then use the abelian localization formula
to give a formula for $F_T (\nusym \pist \si)$ where $\si = e^{\iins \bom}$.
As in (\ref{4.g}),
we define a function $\hh^\nusym(y) $ for $y \in \liets$
by
\beq \label{4.ginj}
\hh^\nusym(y) = \sum_{F \in \calf} \sum_{\a  \in \cala_F}
(-1)^{k_F(\a) } H_{\hat{\g_F}(\a)  } (-y +\mu_T(F) )
\int_F e^\omega \tilde{c}_{F,\a} . \eeq
The notation is as in ({\ref{4.g}) except that
each $\bfjw$ appears in $\hat{\g}_F(\a)$ with multiplicity
$n_\fj(\a)  $ if it is not a root and $n_\fj(\a) - 1 $ if it is a
root. The function
$H_{\hat{\g_F}(\a) } $ is then as defined in Proposition \ref{p3.5}.

In the case
when the action of $T$ has isolated fixed points, the theorem
\ref{t4.1'inj} below
is the main
result
(Theorem 2.16) of Part I  of the paper \cite{GP} of Guillemin
and Prato. For the most part our proof translates directly from
the proof given in that paper, except that we use
the abelian localization theorem in place of
stationary phase.

\begin{theorem} \label{t4.1'inj}
Suppose
that $\lietp$ injects into $T_x X$ for all fixed points
$x$ of the action of $T$.
Then the distribution
$\hh^\nusym$ given in (\ref{4.ginj}) is identical
to  the distribution $ \ggh  = F_T (\nusym \pist e^{\iins \bom}) $.
\end{theorem}
\Proof This theorem is proved in exactly
the same way as   Theorem \ref{t4.1'}.
The conclusion of the proof goes, as in the
proof of Theorem 2.16 of \cite{GP}, by applying
Lemma \ref{l4} directly to
$\hh^\nusym$ as given in
(\ref{4.ginj}) and $\gggh = F_T (\nusym \pist e^{\iins \bom} ) $.
 $\square$

In particular we have the following
\begin{prop} \label{p4.2inj}
Suppose
that $\lietp$ injects into $T_x X$ for all fixed points
$x$ of the action of $T$.
Define
\beq \label{4.3ppinj}
Q(y) = \nusym(y) F_T (\nusym \pist \si) (y), \eeq
so that
\beq \label{4.3pppinj}
\ie = \frac{1}{(2 \pi i )^s  \wn \vol( T) \e^{s/2} }
\tintt [dy] e^{- \inpr{y,y}/{2 \e}  }  Q(y). \eeq
Then $Q$ is a piecewise polynomial function supported on
cones each of which has apex at $\mu_T(F)$ for some component
$F$ of the fixed point set of $T$.
\end{prop}

We now drop the hypothesis that $\lietp$ injects into $T_x X$
at the fixed points $x$ of the action of $T$, and
return to the general situation described in
Proposition \ref{p4.2}.
It will follow from (\ref{5.6}) that $Q$ is smooth near
$y = 0$: thus in particular there is a polynomial $Q_0$ which
is equal to $Q$ near $y = 0$.  Of course $Q_0 = D_\nusym R_0$ where
$R_0$ is the polynomial which is equal to $\ft (\pist \si)$ near
$y = 0$.
In the next section we shall provide an alternative
description of $Q_0$ and prove
\begin{theorem}
\label{t4.3}
$$
\ie_0 \eqdef
\frac{1}{ (2 \pi i)^s   \wn \vol( T) \e^{s/2} }
 \tintt [dy] e^{- \inpr{y,y}/{2 \e} }
\: Q_0(y)  = e^{\e \T} e^{i \om_0} [\xred]. $$
\end{theorem}
This tells us that the contribution to
$\ie$ from $\mu^{-1}(0)/K$ is obtained by integrating
$Q_0$ rather than $Q$ (weighted by the Gaussian $\tilde{g}_{\e^{-1}}$
defined in the first paragraph of this Section)
over $y \in \liets$.

\renorm
\section{The proof of Theorem 4.7}

This section gives the proof of Theorem \ref{t4.3}, which
identifies $\eeth [\xred]$ with the integral of
a Gaussian $\gtsoe$ times a polynomial $Q_0$.
The key step in the proof is
the well known result Proposition  \ref{p5.2} below, which
gives a normal form for the symplectic form, the $K$ action
and the moment map in a neighbourhood $\calo$ of $\zloc$.
We first recast the distribution $Q$ in terms of an integral
over $X$ (Proposition \ref{p5.1}), so that
$Q(y)$ is given by the integral over $X$ of a distribution
supported where $\mu$ takes the value $y$. (This step occurs
also in the  proof \cite{DH} of the Duistermaat-Heckman theorem.)
Hence, for sufficiently small $y$, this distribution
is supported in $\calo$ and we may do the integral over
$X$ to obtain the value of the polynomial $Q_0 $ which
is equal to $Q $ near $0$.
This turns out to be given by an integral over $\xred$ involving
the symplectic form and the curvature of a bundle
over $\xred$ (see (\ref{5.6})).
Finally, we multiply $Q_0$
by the Gaussian $\gtsoe$ and integrate over $\liet$
to see that the result is
 $\eeth [\xred]$.

\begin{prop} \label{p5.1}
$Q(y) = \nusym^2(y) (2 \pi)^{s/2} \int_{x \in X}
e^{i \om} \delta (y - \mu(x) ) $,
where $\delta$ denotes the (Dirac) delta distribution.
\end{prop}
\Proof We have by Lemma \ref{l3.3} that
$$Q = \nusym F_T (\nusym \pist \si ) =
\nusym^2 F_K(\pisk \si) ,$$
so that
$$ Q(y) = \frac{\nusym^2 (y) }{(2 \pi)^{s/2} }
 \intk  [d \phi] \, e^{-i \evab{\phi}{ y} }
\int_{x \in X} e^{i \om} e^{ i \eva{\mu(x)}{ \phi}  } \, $$
 \beq \label{5.1}
= \nusym^2 (y) (2 \pi)^{s/2} \int_X e^{i \om} \delta(\mu - y).  \square \eeq

We would like to study this  for $|y| < h $ for sufficiently
small $h > 0$. Now there is a neighbourhood of $\mu^{-1}(0)$
on which the symplectic form is given in a standard way
related to the symplectic form $\om_0$ on $\xred$: this
follows from the coisotropic embedding theorem
(see   sections 39-41 of \cite{STP}).
\begin{prop} \label{p5.2} {\bf (Gotay\cite{gotay}, Guillemin-Sternberg
\cite{STP}, Marle \cite{marle})}
Assume $0$ is a regular value of $\mu$ (so that $\mu^{-1}(0)  $
is a smooth manifold and  $K$ acts on $\mu^{-1}(0) $
with finite stabilizers).
Then there is a neighbourhood $\calo \cong \mu^{-1}(0) \times
\{ z \in \lieks, |z| \le h \} $
$\subseteq \mu^{-1}(0) \times \lieks$  of $\mu^{-1} (0)$ on which the
symplectic form is given as follows. Let
$P \eqdef \mu^{-1}(0) \stackrel{q}{\to} \xred $
be the orbifold  principal $K$-bundle given by
the projection map $q: \zloc \to \zloc/K$,
 and let $\theta$ $ \in \Om^1(P) \otimes \liek$
be a connection
for it. Let $\omr$   denote the induced symplectic form on
$\xred$, in other
words $q^* \om_0 = i_0^* \om$.
Then if we define a 1-form
$\tau$ on $\calo\subset P \times \lieks$ by
$\tau_{p,z}  = z(\theta)$ (for $p \in P$ and $z \in \lieks$),
  the symplectic form on $\calo$
is given by
\beq \label{5.2} \om = q^* \omr + d \tau.  \eeq
Further,  the moment map on $\calo $
is given by $\mu (p, z) = z$.
\end{prop}

\noindent{\em Proof of Theorem \ref{t4.3}:}~
We assume for simplicity of notation that $K$
acts freely on $\mu^{-1}(0)$, but all of the following
may be transferred to the case when $K$ acts with finite
stabilizers by introducing $V$-manifolds or orbifolds (see \cite{kaw}).
In other words, we work locally on  finite covers
of subsets of $\mu^{-1}(0)$ and $\mu^{-1}(0)/K$, where the
covering group is the stabilizer of the $K$ action at a
point $x \in \mu^{-1}(0)$.

When  $|y| < h$ and $h$ is sufficiently small,
the distribution $\delta(\mu(x)  - y)$
is supported in $\calo$, so we may compute
$Q_0(y)$ from (\ref{5.1}) by restricting to $\calo$.
We have
\beq \label{5.3}
Q_0(y) = (2 \pi)^{s/2}  \nusym^2(y)  \int_{(p, z') \in P \times \lieks}
e^{i \om} \delta(y - z')  \eeq
\beq \label{5.4}
 \onebl = (2 \pi)^{s/2} \nusym^2(y)
\int_{(p, z') \in P \times \lieks}
\exp i(q^* \om_0 + \evab{d \t}{ z'})  \exp  i\evab{\t}{ dz'} \delta(y - z')
\eeq
Now the term in  $\exp  i \evab{\t}{ dz'} $ which contributes
to the integral (\ref{5.4}) is $i^s \Om \,[dz'] $
where $[dz'] $ is the volume form on $\liek$
(since all factors $d z'_1 \dots {dz'}_l$ must appear in
order to get a contribution to the integral).
Here, $\Om  = \prod_{j = 1}^s \t^j$ (for $j$ indexing an
orthonormal basis of $\liek$ and $\t^j$ the
corresponding components of the connection
$\t$) is a form integrating to $\vol(K)$
over each fibre of $P \to \xred$.

Doing the integral over $z' \in \lieks$, we get
\beq \label{5.5}
Q_0(y) = i^s \nusym^2(y) (2 \pi)^{s/2}
\: \int_P \exp	i (q^* \omr + \evab{d \t +
[\t, \t]/2 }{ y }
 ) \: \Om \eeq
\beq \label{5.6}
\bla = i^s  \nusym^2 (y) (2 \pi)^{s/2}  \: \int_{\zloc} \exp
i
(q^* \omr +  \evab{F_\t }{ y } ) \, \Om.
\eeq
Here, $F_\t = d \t + \half[\t, \t]$ is the curvature
associated to the connection $\t$;
we may introduce the term $[\t, \t]$ into the
exponential in (\ref{5.5}) since
the additional factors $\t$  will give zero under
the wedge product with $\Om$. Formula (\ref{5.6})
shows that $Q_0$ is a polynomial in $y$.

Now we were interested in
$$\ie_0 = \frac{1}{(2 \pi i )^s   \wn \vol( T) \e^{s/2} }
\tintt [dy] e^{- \inpr{y,y}/{2 \e} } Q_0(y) $$
\beq \label{5.7}
= \frac{1 }
{ (2 \pi)^{s/2}  \wn   \vol( T) \e^{s/2} }
\tintt [dy] \nusym^2 (y) \:  e^{- \inpr{y,y}/{2 \e} } \:
\int_{P} \exp (q^* \omr + \evab{F_\t}{ y} ) \, \Om, \eeq
\beq \label{5.9}
= \frac{1}{ (2 \pi \e)^{s/2}  \vol(K)  } \int_{z \in \lieks}
\, [dz] \, e^{- \inpr{z,z}/{2\e}  } \: \int_{P}
\exp (q^* \omr + \evab{F_\t}{z} ) \, \Om, \eeq
where the last step uses Lemma \ref{l3.1} and the
fact that $\int_{P} \exp (q^* \omr + \evab{F_\t}{ z} )  \Om$
is an invariant function of $z$.

We now regard $F_\t$ as a formal parameter and complete the
square to do the integral over $z$: we have
(identifying $z(F_\t)$ with $\inpr{F_\t, z}$  using the invariant
inner product $\inpr{\cdot, \cdot}$)
\beq \label{5.10}
\int_{z \in \lieks} [dz] e^{- \inpr{z,z}/{2 \e} }
\: \exp \inpr{F_\t, z}
 =  (2 \pi \e)^{s/2} \exp \e \inpr{F_\t, F_\t}/2. \eeq
But $\inpr{F_\t, F_\t}/2$ is just the class
$\pi^* \T$ on $P$, for $\T \in H^4 (\xred)$.
Hence we obtain (integrating over the fibre of
$P \to \xred$ and using the fact that the
integral of $\Om$ over the fibre is $\vol(K)$)
\beq \label{5.11}
\ie_0 = \int_{\xred} \exp i\om_0 \: \exp \e \T, \eeq
completing the proof of Theorem \ref{t4.3}. $\square$

\renorm
\section{The proof of  Theorem 4.1 }

In this section we complete the proof of Theorem
\ref{t4.1}. This is done by observing
that $\ie - \eeth[\xred]$ is of the
form $\ims \int_{\liet} \gtsoe (Q - Q_0) = $
$\int_{\liet} \gtsoe D_\nusym ( R - R_0)$ =
$\int_{\liet} (D_\nusym^* \gtsoe) (R-R_0)$, where
$R - R_0$ is piecewise polynomial and
supported away from $0$. (Here,
$D_\nusym^* = (-1)^{(s-l)/2} D_\nusym.$)
The results of  \cite{Ki1}
establish that the   distance
of any point of Supp($Q - Q_0$) from $0$
is at least $|\b|$ for some nonzero $\b$
in the indexing set $\calb$ defined in Section 4.
Hence we obtain the estimates in Theorem \ref{t4.1}.

In fact the function $R - R_0$ is known explicitly
in terms of the values of $\mu_T(F)$ (where $F$ are
the components of the fixed point set), the
integrals over $F$ of characteristic classes of subbundles
of the normal bundle $\nu_F$, and the weights of
the action of $T$ on $\n_F$ (see (\ref{2.4}) and
Proposition\ \ref{p3.5}). The function $R - R_0$ is
polynomial on polyhedral regions of $\liet$,
so that the quantity  $\ie -  \eeth [\xred]$ can in principle
be computed from the integral of a polynomial times a Gaussian
over these polyhedral regions.
We shall study these integrals in another paper and
relate them to the cohomology of the higher strata
in the stratification of $X$ according to the gradient
flow of $|\mu|^2$ given in \cite{Ki1}.

We now   examine $\ie - \ie_0$ and prove Theorem \ref{t4.1}.
Recall from section 4 that the indexing
set $\calb$ of the critical
sets $C_\b$ for the
function $\rho = |\mu|^2$ is
 $\calb = \liet_+ \cap W\calb$ where
$W \calb = \{ w \b: \; \b \in \calb, \onebl
w \in W \}$ is the
set of closest points to $0$
of convex hulls of nonempty subsets
of the set $\{ \mu_T (F): F \in \calf\}$
of images under $\mu_T$ of the connected
components of the fixed point set of $T$ in $X$.
We shall refer to $\{ \mu_T (F): F \in \calf\}$  as the
set of {\em weights}\footnote{The motivation for this is
that if $X$ is a nonsingular subvariety of complex
projective space ${ \Bbb P}_n$ and $T$ acts
on $X$ via a linear action on $\CC^{n+1} $
then each $\mu_T(F)$ is a weight of this action
when appropriate identifications are made.}  associated to $X$ equipped with
the
action of $T$.

Let $\calj$ denote the locus
$$\calj = \{ y \in \liets  : \, Q \onebl
\mbox{is not smooth at} \onebl y\}. $$
Then we have
\begin{prop}
$\calj \subset \caljab, $ where
$\caljab = \{ y \in \liets : \, \ft
(\pis e^{\iins \bom})  $ is not smooth at $y \}. $
\end{prop}
\Proof $Q = \nusym \ft (\nusym \pis e^{\iins \bom})$
$ = \nusym D_\nusym
\ft(\pis e^{\iins \bom})$.
Hence if $\ft (\pis e^{\iins \bom})$ is smooth at $y$ then so is
$Q$. $\square$

Now it follows from \cite{JGP} (Section 5) that
$\caljab = \cup_{\g \in \G} \mu_T (V_\g)$ where
$V_\g$ is a component of the fixed point set of a one parameter subgroup
$T_\g$ of $T$ and $\G$ indexes all such one parameter subgroups and
components of their fixed point sets.
Let $$D = \cap \{ D_\b: \;  \b \in W\calb - \{0\}
   \}$$
where $D_\b$ denotes the open half-space
$$ D_\b = \{ y \in \liets: \;  y(\beta)
<  |\b|^2 \}. $$
Note that if $\b \in \calb -  \{0\}$
then $D_\b$ contains $0$ and its boundary
is the hyperplane
$$H_\b = \{ y \in \liets: \; y(\beta) = |\b|^2 \}.$$
\begin{lemma} \label{l6.2}The support of $Q - Q_0$ is contained in
the complement of $D$ (or equivalently
$Q = Q_0$ on $D$).
\end{lemma}
\Proof
Suppose $V_\g$ is a component of the fixed point set of a one
parameter subgroup $T_\g$.
By the Atiyah-Guillemin-Sternberg convexity theorem
\cite{aam,gsconv},
$\mu_T(V_\g)$ is the convex hull of some subset of the
weights. Hence the closest point to $0$ in
$\mu_T(V_\g)$ is in $W \calb$.
Now either this closest point is $0$, 
or else the closest point is a point
$\b \in
W\calb - \{0\}$ (in which case $\mu_T(V_\g) \subset \liet - D_\b
\subset \liet - D$).

Now if $x$ is a point in Supp($Q - Q_0$) then the ray from
$0$ to $x$ must pass through at least one point in $\calj$: hence
it suffices to prove $\calj \subset \liet - D$ since $D$ is
the intersection of a number of
open half spaces all of which contain $0$.

But $\calj \subset \cup_\g \mu_T (V_\g)$ and every
point in $\mu_T(V_\g)$ is in $\liet - D$ unless $0$ is
in $\mu_T(V_\g)$. Moreover if $x \in \mu_T(V_\g) \cap \calj$
and $0 \in \mu_T(V_\g)$, we may consider the function
$Q - Q_0$ restricted to a small neighbourhood of the ray
from $0$ through $x$. This ray lies in the hyperplane
$\tilde{H}_\g$ which
is the orthocomplement in
$\liet$ to the Lie algebra $\liet_\g$
of $T_\g$. (Since the component of
$\mu_T$ in the direction of $\liet_\g$ is
constant along $V_\g$ and since
$0 \in \mu_T(V_\g)$, it follows
that $\mu_T(V_\g)$ is contained in
$\tilde{H}_\g$.)  Since
$Q - Q_0$ is identically zero near $0$ but not near $x$,
the ray from $0$ to $x$ must contain a point $x'$ in
$\calj$ which is contained in some $\mu_T(V_{\g'})$
with $t_\g \ne t_{\g'}. $ If $0 \notin
\mu_T(V_{\g'})$
  then $\mu_T(V_{\g'} ) \subset \liet - D$
and so $x \in \liet - D$. If $0 \in
\mu_T(V_{\g'}) $ then we simply
repeat the argument, considering the restriction
of $Q - Q_0$ to a neigbourhood
of $\tilde{H}_\g \cap \tilde{H}_{\g'}$.
Since $0 \notin \calj$, after finitely many repetitions of
this argument we get the required result.
 Hence the Lemma is proved. $\square$

To complete the proof of Theorem \ref{t4.3},
we then use Lemma \ref{l6.2} to express $\ie - \ie_0 $ as
\beq \label{6.3}
\ie - \ie_0
 =  \frac{1}{(2 \pi i)^s |W|   \vol T \e^{s/2} }
\int_{\liets - D} [dy] (Q - Q_0) e^{- |y|^2/{2 \e} } . \eeq
Denote by $C$ the set
$\{ y \in \liets  - D  \: : \: |Q(y) - Q_0(y) | \le 1 \}. $
Then $$(2 \pi)^s |W| \vol(T) \e^{s/2}
|\ie - \ie_0| \le  \int_C [dy]\stg +
\int_{\liets - D}[dy]  |Q - Q_0|^2 \stg. $$
If $b $ is the minimum value of $|\beta | $
over all $\beta \in \calb - \{0\}$, then
$$\int_C [dy]\stg  \le \int_{|y| \ge b } [dy] \stg
\le e^{- b^2/{2 \e} } q(\e), $$
where $q(\e) $ is a polynomial in $\e^\half$.

Further, denote by $p$ the function  $|Q - Q_0|^2$.
Then
\beq \label{6.4}\int_{\liets - D} [dy] p(y) \stg
  \le \sum_{\beta \in W\calb - \{0\}}
\int_{y \in D_\beta}[dy] \stg p(y). \eeq
For each $\b \in W \calb - \{0\} $ one can now decompose
$y \in \liets$ into $y = w_0 \hatb + w, w \in \beta^\perp$, $w_0 \in \RR$
(where $\hatb = \b/|\b|$).
Hence each of
the integrals (\ref{6.4})  is of the form
$$ \int_{w_0 \ge |\beta| } \int_{w \in \beta^\perp}
e^{- w_0^2/{2 \e} } e^{ - |w|^2/{2 \e} } p(w_0, w), $$
and this is clearly bounded by $e^{- |\beta|^2/{2 \e} } $ times
a polynomial in $\e$. This completes the proof of Theorem 4.1.
$\square$

\renorm
\section{Extension of Theorems 4.1 and 4.7 to other classes}

In this section we extend Theorems \ref{t4.1} and \ref{t4.3}
to equivariant cohomology classes of the form $\z = \eta e^{\iins \bom} $
where $\eta \in \hk(X)$.
 More precisely
we shall show the following
\begin{theorem} \label{t7.1}
 Suppose $\eta \in \hk(X)$
and suppose
that
$i_0^* \eta \in \hk (\zloc)$ is represented by $\eta_0 \in H^*(\zloc/K)$
(where $i_0: \zloc \to X$ is the inclusion map). Then

{\bf(a)} We have that
$$ \eta_0 e^{\e \T}  e^{i \om_0}
[\xred]  = \frac{1}{(2 \pi i )^s |W| \vol(T) \e^{s/2} }
\tintt [dy] e^{- |y|^2/{2 \e} } \,  Q_0^\eta (y), $$
where $Q^\eta(y) = \nusym(y) \ft (\nusym \pist (\eta \exp \iins \bom) ) $
and $Q_0^\eta(y) $ is a polynomial
which is equal to $Q^\eta(y) $ near $y = 0$.

{\bf (b)}
Let
$ \rho_\beta = |\beta|^2$ be  the
value of the function $|\mu|^2: X \to \RR$ on the
critical set $C_\b$.   Then
there exist functions $h_\beta: \RR^+ \to \RR $  such that
for some $N_\beta \ge 0$,
$\e^{N_\beta} h_\beta(\e)$ remains bounded as $\e \to 0^+$, and
for which
$$\Bigl | \,  \frac{1}{(2 \pi i )^s \vol(K) } \intk  [d \phi]
e^{-\e |\phi|^2/2} \pisk \eta e^{\iins \bom}
 - \eta_0 e^{\e \T} e^{i \om_0}  [\xred]  \,  \Bigr | \le
\sum_{\beta \in \calb - \{0\} }
 e^{- \rho_\beta  /{2 \e}  } \, h_\beta(\e). $$

{\bf (c)} Suppose $\eta$ is represented by
$\tilde{\eta} = \sum_J \eta_J \phi^J $ $ \in \Om^*_K(X)$
for
$\eta_J \in \Om^*(X)$. Then $Q^\eta$ is of the form
$Q^\eta (y) = \sum_J D_J R_J(y)$, where $D_J$ are
differential operators on $\liets$ and $R_J$ are piecewise
polynomial functions on $\liets$.
\end{theorem}
\noindent{\em Proof of } {\bf(a)}: We examine the function $ \qeta(y) $ near
$y = 0$. We have from Lemma \ref{l3.3} and the paragraph before
Theorem \ref{t2.1} that
$$ \nusym \ft (\nusym \pist (\eta e^{\iins \bom} ) ) (y)
 = \Bigl (\nusym^2 \fk  \pisk (\eta e^{\iins \bom} ) \Bigr )  (y)$$
$$ = \frac{\nusym^2 (y)}{(2 \pi)^{s/2} }  \intk [d \phi] \int_{x \in X }
e^{- i \eva{y}{ \phi} } e^{i \om} e^{ i \eva{\mu(x)}{ \phi} } \eta (\phi). $$
Now since $\eta (\phi) = \sum_I \phi^I \eta_I$ for $
\eta_I \in \Om^*(X)$  (where the $I$ are multi-indices), we may
define for  any $x \in X$ a distribution $\cals_x$ with
values in $\L^* T^*_x X$ as follows: for any $y \in \liets$ we have
$$\cals_x(y) = \intk [d \phi] \eminev{y}{ \phi }   e^{i \om}
\epinev{\mu(x)}{ \phi} \sum_I \eta_I \phi^I $$
\beq  = \sum_I (i \partial/\partial y)^I \intk [d \phi]
\epinev{ (\mu(x) - y)}{ \phi} e^{i \om} \eta_I \eeq
$$ = (2 \pi)^s
\sum_I (i \partial/\partial y)^I \delta(\mu(x) - y) e^{i \om} \eta_I. $$
Thus the distribution $\cals_x(y)$, viewed as
a distribution $  \cals (x,y)   $ on $X \times \liets$,
 is supported on $\{ (x , y) \in X \times \liets   \, | \,
\mu(x) = y \}. $ Hence for sufficiently small $y$,
$\cals(x,y) $ (viewed now
as a function of $x$)   is supported on $ \calo$ (in the notation
of Section 5) and we find that
\beq \qeta(y)  = \frac{\nusym^2(y)}{\htps}
\int_{x \in X} \cals(x, y) \eeq
$$ = \frac{\nusym^2(y)}{ (2 \pi)^{s/2} }  \intk [d \phi]  \int_{P \times
\lieks}
\epinev{(\mu(x)  - y)}{ \phi}  e^{i \om} \eta(\phi)  . $$

Now consider the restriction of $\eta$ to $\hk(P \times \lieks) \cong
\hk(P)$ (where $P = \mu^{-1}(0)$ as in Section
5). Recall that  there exists $\eta_0 \in
\Om^*(P/K)$ such that $\eta - \pi^* q^*\eta_0 = D \gamma$ for some
$\gamma$, where $D$ is the equivariant cohomology differential on
$P \times \lieks$ and $\pi: P \times \lieks \to P  \times \{0\}$ and
 $q: P \to \to P/K$ are the
projection maps.\footnote{This is because the map
$i \circ \pi: P \times \lieks \to P \times \lieks$ is homotopic
to the identity by a homotopy through equivariant maps, where
 $i: P \times \{ 0 \} \to
P \times \lieks$ is the inclusion map. Hence
$i$ induces an isomorphism $i^*: \hk (P \times \lieks)
\to \hk(P \times \{0\}). $}
 We then have that
\beq \label{7.d3}
Q^\eta (y) - \frac{\nusym^2(y)}{(2 \pi)^{s/2} } \intk [d \phi]
\int_{x \in
P \times \lieks} \epinev{(\mu(x) - y)}{\phi} e^{i \om} \pi^* q^*\eta_0  \eeq
$$ \onebl \onebl = \frac{\nusym^2(y)}{(2 \pi)^{s/2} } \intk [d \phi]
\int_{x \in P \times \lieks} \epinev{(\mu(x) - y)}{ \phi} e^{i \om}
D \gamma \eqdef \frac{\nusym^2(y)}{(2 \pi)^{s/2}}  \bigtriangleup. $$
But also
$$ \epinev{(\mu(x) - y)}{ \phi} e^{i \om} D \gamma = D (
\epinev{(\mu(x) - y)}{ \phi}
e^{i \om} \gamma ) $$
(since $D \phi_j = 0 $  and $D \bom = 0$). Hence
$\bigtriangleup = \intk [ d \phi] \int_{P \times \lieks}
d (\epinev{(\mu(x) - y)}{ \phi} e^{i \om} \gamma) $ (since for any
differential form $\Psi$, the term
$\iota_{\tilde{\phi}}  \Psi$ in $D \Psi$
cannot contain differential forms of top degree in $x$).
Using Stokes' Theorem and replacing $P \times \lieks$ by
$P \times B(\lieks) $ where $B(\lieks)$ is a large ball in
$\lieks$ with boundary $S(\lieks)$, we have that
$$ \bigtriangleup =
\intk [d \phi] \int_{P \times S(\lieks) }
\epinev{(\mu(x) - y)}{ \phi} e^{i \om} \sum_I \gamma_I \phi^I$$
for $\gamma_I \in \Om^*(P \times \lieks)$. Thus we have
$ \bigtriangleup = \sum_I  (i \partial/\partial y)^I S_I(y) $ where
$$ S_I (y) =
\intk [d \phi] \int_{x \in P \times S(\lieks) }
\epinev{(\mu(x) - y)}{ \phi} e^{i \om} \gamma_I. $$
Now we do the integral
over $\phi$ to obtain
$$ S_I(y)  = (2 \pi)^{s}  \int_{x \in P \times S(\lieks)} \delta (\mu(x) - y)
e^{i \om} \gamma_I. $$
This is zero since the delta distribution is supported
off $S(\lieks)$ (recall that $\mu(p, z) = z$ for
$(p, z) \in P \times \lieks$). Hence we have that
$\bigtriangleup = 0$, and so
\beq  \qeta(y) = \frac{\nusym^2(y)}{(2 \pi)^{s/2} } \intk [d \phi]
\int_{(p, z) \in P \times \lieks} \epinev{(z - y)}{ \phi}
e^{i \om} \pi^* \eta_0, \eeq
and the rest of the proof is exactly the same as the proof of
Theorem \ref{t4.3} which was for the case $\eta_0 = 1$.
In particular the analogue of (\ref{5.6}) is
\beq \label{5.6ext}
Q^\eta_0(y)  = i^s \nusym^2 (y) (2 \pi)^{s/2}  \: \int_{\zloc}
\eta_0 \exp
(q^*\omr +  \evab{F_\t}{  y } ) \, \Om \, ;
\eeq this equation
shows that $Q_0^\eta$ is a polynomial (and in
particular smooth) in $y$.

\noindent{\em Proof of }{\bf(b):}   This  is a direct extension of the
proof of Theorem \ref{t4.1}, with $Q$ and $Q_0$ replaced by
$\qeta$ and $Q^\eta_0$.

\noindent{\em Proof of }{\bf (c):}  Since $\eta = \sum_J \eta_J \phi^J$,
the  abelian localization formula for
$\pis (\eta e^{\iins \bom})$ yields
\beq \label{7.6n}
\nusym(\phi) \pis (\eta e^{\iins \bom}) (\phi) =  \nusym(\phi)
\sum_{F \in \calf} \rfe(\phi), \bla
\rfe(\phi) =  \int_F \frac{i_F^* \eta( \phi) e^{i \bom(\phi)} }{e_F(\phi)},
\eeq $$ \bla =
\sum_J \phi^J
e^ {i \eva{\mu_T(F)}{ \phi} }
 \int_F \frac{i_F^* \eta_J e^{i \om}
  }{e_F(\mar \phi) } $$
$$ \onebl = \nusym(\phi) \sum_J \phi^J
\sum_{F \in \calf, \a \in \cala_F} \tfaj(\phi) , $$
\beq \label{7.06}\mbox{where} \onebl \onebl  \tfaj(\phi) =
 (-1)^{k_F(\a) }
\int_F i_F^* \eta_J
e^{i \om}
\tilde{c}_{F,\a}
   \frac{ \epinev{ \mu_T(F)}{ \phi }   }
{\prod_j \Bigl ( \bfjw (\mar \phi) \Bigr )^{{n}_\fj(\a)} }
  \eeq
(and the $\tilde{c}_{F,\a}$ are as in
(\ref{4.004})).
Here, as in Section 4, we may form a distribution
\beq \label{7.07}H(y) = D_\nusym   \tilde{H},
\onebl \tilde{H} = \sum_J (i \partial/\partial y)^J
\sum_{F \in \calf, \a \in \cala_F} (-1)^{k_F(\a)}
\Bigl ( \int_F  i_F^* \eta_J
e^{i \om}
\tilde{c}_{F,\a} \Bigr )  H_{\bar{\g}_F(\a) } ( \mu_T(F) - y ), \eeq
where the piecewise polynomial function
$H_{\bar{\b} }$ is as in Proposition \ref{p3.5}. Thus $\tilde{H}$ is of
the form $\tilde{H}(y) = \sum_J D_J R_J (y)$ where
the $R_J(y)$ are piecewise polynomial and
the $D_J$ are differential operators.

We also define the distribution
\beq \label{7.07'}
G(y) =  \ft( \nusym \pis \eta e^{\iins \bom}) . \eeq
Thus $\nusym(y) G(y) = Q^\eta (y)$.

Then we may show
that  the distributions $G$ and $H$ are identical.
For we apply Lemma \ref{l4} as before. We find by Proposition
\ref{p3.5}(a) that $(\ft G)(\psi) $ and $(\ft H)(\psi)  $ are identical
off the hyperplanes $\bfj( \psi) = 0$, so the first
hypothesis of Lemma \ref{l4} is satisfied. Moreover,
$H$ is supported in a half space and
\beq G(y) = \nusym(y) \fk (\pis \sum_J \eta_J \phi^J e^{\iins \bom} ) (y) \eeq
$$ = (2 \pi)^s \nusym(y) \sum_J (i \partial/\partial y)^J \int_{x \in X}
\eta_J  \delta(\mu(x) - y) e^{i \om}, $$
so $G$ is supported in the compact set $\mu(X) \cap \liets$ and hence
$H - G$ is supported in a half space. Thus the second
hypothesis of Lemma \ref{l4} is also satisfied and
we may conclude that $G = H$, which completes the proof.
$\square$

\renorm
\section{Relations in the cohomology ring of symplectic quotients}
In this section we shall prove a  formula
for the evaluation of cohomology
classes from $\hk(X) $ on the fundamental
class of $\xred$, and apply
it to study the cohomology
ring $H^*(\xred)$ in two examples.

\nc{\gamp}{\G(P) }
\nc{\npl}{ {n_+} }

We shall prove the following theorem:

\begin{theorem}\label{t8.1}
Let $\eta \in \hk(X)$ induce $\eta_0 \in H^*(\xred)$. Then we have
\beq \label{8.00}
\eta_0  e^{i \om_0} [\xred]  =
\frac{(-1)^\npl }{ (2 \pi)^{s-l} |W| \vol(T) }
\treso
\Biggl (     \: \nusym^2 (\psi)
 \sum_{F \in \calf} \rfe(\psi) [d \psi] \Biggr ), \eeq
where
$$ \onebl \rfe(\psi) =  e^{i \mu_T(F) (\psi) }
\int_F \frac{i_F^* (\eta(\psi) e^{i \omega})  }{e_F(\mar  \psi) }
 . $$
\end{theorem}
Here $s$ and  $l$ are the dimensions of $K$ and
its  maximal torus  $T$,
and $\nusym(\psi) = \prod_{\g > 0 }
\g(\psi)$ where $\g$ runs over the positive roots of $K$;
the number of positive roots
$(s - l)/2$  is denoted by $\npl$.
If $F \in \calf$ (where $\calf$ denotes the set
of components of the $T$ fixed point set) then
$i_F: F \to X$ is the inclusion of $F$ in $X$ and
$e_F $ is the equivariant Euler class of the
normal bundle to $F$ in $X$.
The quantity $\treso (\Omega)$
will be defined
below (Definition \ref{d8.5n}).
The definition  of $\treso (\Omega)$
 will depend on the choice of a cone $\lasub$, a test function
$\testf$, and
 a ray in $\liets$ specified by
a parameter  $ \zray \in \liets$ -- but in fact if the form $\Omega$
is sufficiently well behaved
(as is the case in (\ref{8.00}))  then the
quantity $\treso (\Omega) $ will turn out to be
independent of these choices. (See  Propositions
\ref{p8.6n}, \ref{p8.7n} and \ref{p8.8n}.)

In the case of $K = SU(2) $  the result of
Theorem  \ref{t8.1} is as follows:

\begin{corollary} \label{c8.2}
Let $K = SU(2)$, and let $\eta \in \hk(X)$ induce
$\eta_0 \in H^*(\xred)$.
Then
the cohomology class $\eta_0   e^{i \om_0}  $ evaluated on
the fundamental class of $\xred$ is given by the
following formula:
$$ \eta_0  e^{i \om_0}  [\xred] =
-\frac{ 1}{2} {\rm Res}_0
\Biggl (  \psi^2
\sum_{F  \in \calf_+ } \rfe(\psi) \Biggr ), \onebl
{\rm where} \onebl   \rfe (\psi) =
e^{i \mu_T(F)(\psi) }
\int_F
\frac{i_F^* \eta(\psi) e^{i \omega}   }{e_F( \psi) }.
$$ Here,
${\rm Res_0}$ denotes the coefficient of $1/\psi$, and
$\calf_+$ is the subset of the fixed point set of $T = U(1)$
consisting of those components $F$ of the $T$ fixed point set
for which $\mu_T(F) > 0 $.
\end{corollary}
An important  special case (cf. Witten \cite{tdg}, Section 2.4) is as follows:
\begin{corollary} \label{c8.00}
Let $\eta \in \hk(X)$ induce $\eta_0 \in H^*(\xred)$,
and let $\Theta \in H^4(\xred)$ be induced by the
polynomial function $- \inpr{\phi, \phi}/2$
of $\phi$, regarded as an element of $H^4_K$.
Then if $\e > 0$ we have
\beq \label{8.010}
\eta_0 e^{\e \Theta} e^{i \om_0} [\xred]  =
 \frac{(-1)^\npl }{(2 \pi)^{(s-l)} |W| \vol(T) }
\sum_{m \ge 0} \frac{1}{m!} \treso
\Biggl (    (- \e |\psi|^2/2)^m  \: \nusym^2 (\psi)
 \sum_{F \in \calf}  \rfe(\psi) [d \psi] \Biggr ), \eeq
$$  \bla {\rm where} \bla  \rfe (\psi) = e^{i \mu_T(F) (\psi) }
\int_F \frac{i_F^* (\eta(\psi) e^{i \om})  }{e_F(\mar  \psi) }
. $$
\end{corollary}

We  now give a general definition (Definition
\ref{d8.5n}) of $\ress h [d \psi]$
when $h $ is a meromorphic function on an open
subset of $\liet \otimes \CC$ satisfying
certain growth conditions at infinity. For the restricted class of
meromorphic forms of the form
$e^{i \l(\psi)} [d \psi]/\prod_j \b_j (\psi), $ Definition \ref{d8.5n}
implies
the existence of an explicit procedure for computing these residues
by successive contour integrations, which is outlined in
Proposition \ref{p8.4}.

Our definition of the residue will be based on the following
well-known result (\cite{hor}, Theorem 7.4.2 and Remark following
Theorem 7.4.3):
\begin{prop} \label{cty}
(i) Suppose $u$ is a distribution on $\liets$. Then
the set $\gu = \{ \xi \in \liet: \: e^{(\cdot,  \xi) } u $
is a tempered distribution$\}$ is convex. (Here, $(\cdot, \cdot)$
denotes the pairing between $\liet$ and $\liets$.)

\noindent (ii) If the interior $\guo$ of $\gu$ is nonempty,
then there is an analytic function $\hat{u} $ in
$\liet + i \guo$ such that the Fourier transform
of $e^{(\cdot,  \xi) } u $ is $\hat{u} (\cdot + i \xi)$
for all $\xi \in \guo$.

\noindent (iii) For every compact subset $M$ of $\guo$
there is an estimate
\beq  \label{bds} | \hat{u} (\z) | \le C (1 + |\z|)^N, \onebl
{\rm Im} (\z) \in M. \eeq

\noindent (iv) Conversely if $\G$ is an open convex set
in $\liet$ and $h$ is a holomorphic function on $\liet + i \G$
with bounds of the form (\ref{bds}) for every compact
$M \subset \G$, then there is a distribution
$u$ on $\liets$  such that $e^{(\cdot, \xi)} u $ is a tempered distribution
and has Fourier transform
$h (\cdot + i \xi)$ for all $\xi \in \G$.

\noindent (v)  Finally, if $u$ itself is tempered then the Fourier transform
$\hat{u}$
 is the limit  (in the space $\cals'$
of tempered distributions)
of the distribution  $ \psi \mapsto
\hat{u} (\psi + i t \theta)$
as $t \to 0^+$, for any $\theta \in \guo$.
\end{prop}

\begin{definition} \label{d8.5n}
Let $\lasub$ be a (proper) cone in $\liet$. Let $h$ be a holomorphic
function on $\liet - i {\rm Int} (\L) \subseteq \liet \otimes \CC $
  such that
for any compact subset $M$ of $\liet - i {\rm Int} (\lasub) $ there
is   an integer $N \ge 0 $ and a  constant $C$ such that
$|h (\z)| \le C (1 + |\z|)^N$
 for  all   $\z \in M$.
Let $\testf: \lieks \to \RR$ be a smooth invariant function with
compact support and strictly positive  in some  neighbourhood of $0$,
 and let $\htestf =
\fk \testf: $ $\liek
\to \CC$ be its Fourier transform. Let $\htestfe(\phi) =
\htestf (\e \phi)$ so that $(\fk \htestfe)(z) = \testfe (z) =
\e^{-s} \testf (z/\e). $ Assume $\htestf(0) = 1. $
Then we define
\beq \label{8.1n} \reslch  (h [d \psi])
 = \lim_{\e \to 0^+} \frac{1}{(2 \pi i )^l }
\int_{\psi \in \liet - i \xi} \htestf (\e \psi) h (\psi) [ d\psi] \eeq
where $\xi$ is any element of $\L$.
 \end{definition}

By the Paley-Wiener theorem (\cite{hor}, Theorem 7.3.1), for any
fixed $\xi \in \liek$ the function $\htestf (\psi - i \xi) =
\fk(e^{- \inpr{\xi, \cdot} } \testf)(\psi) $ is rapidly decreasing
since $\fk \htestf = \testf$ is smooth and compactly supported. Hence
the integral (\ref{8.1n}) converges. Now the function $\htestf$ extends
to a holomorphic function on $\liek \otimes \CC$ and in
particular on $\liet \otimes \CC$ (Proposition \ref{cty}), and
by assumption  $h$ extends
to a holomorphic function on
$\liet - i \, {\rm Int} (\lasub)$;
 hence, by Cauchy's theorem, the integral
(\ref{8.1n}) is  independent of $\xi \in {\rm Int}  (\lasub)$.


The independence of $\reslch (\Omega)$ of the
choices of $\lasub$ and $\testf$ when $\Omega$ is sufficiently
well behaved is established by the next results.
\begin{prop} \label{p8.6n}
 Let $h: \liek \to \CC$ be a
$K$-invariant function. Assume that
$\fk h$ is compactly supported; it then follows that $h: \liek \to \CC$
is smooth (\cite{hor}, Lemma 7.1.3) and extends to a
holomorphic function on $\liek \otimes \CC$ (Proposition \ref{cty}).
 Then $\reslch ( h [d \psi]  ) $
is independent of the cone $\lasub$.
\end{prop}

\Proof As above, define $\testfe (z) = \e^{-s} \testf (z/\e),$
so that $\htestfe (\phi) = \htestf(\e \phi). $
The function $h$ extends in particular to a holomorphic
function on $\liet \otimes \CC$, and by the remarks after
Definition \ref{d8.5n}, the function $\htestfe$ also extends to a
holomorphic function on $\liet \otimes \CC$.
Hence Cauchy's theorem shows that
for any choice of the cone $\lasub$,
\beq \label{8.4n}
\reslch ( h(\psi) \dps ) =  \lim_{\e \to 0^+}
\frac{1}{(2 \pi i)^l}
 \intt \htestf (\e\psi) h(\psi) [d \psi]. \bla \square    \eeq
\noindent{\em Remark:} If $h$ is as in the statement of
Proposition \ref{p8.6n},  then $\nusym^2 h$ also satisfies the hypotheses
of the Proposition, so
$\reslch(\nusym^2 h [d \psi] ) $ is also independent of the cone
$\lasub$.

The following is a consequence of Proposition \ref{cty}:

\begin{prop} \label{p8.7n}
 Let $u: \liets \to \CC$ be a distribution,
and assume the set $\gu$ defined in Proposition \ref{cty} contains
$- {\rm Int} (\L)$. Thus $h = \ft u$ is a holomorphic
function on $\liet - i {\rm Int} \L$, satisfying the
hypotheses in Definition \ref{d8.5n}. Suppose in addition
that $u$ is smooth at $0$. Then
$\reslch (h \dps) $ is independent of the test function $\chi$, and
equals $i^{-l} u(0)/(2 \pi)^{l/2} $.
\end{prop}

We shall be dealing with functions $h: \liek \to \CC$ whose Fourier
transforms are smooth at $0$ but which are sums of other
functions not all of whose Fourier transforms need be smooth
at $0$. For this reason we must introduce a small generic parameter
$\zray \in \liets$ where all the functions in this sum are smooth.
More precisely we make the following
\begin{definition} \label{d8.7n} Let $\lasub, $ $\testf$ and $h$ be as in
Definition \ref{d8.5n}. Let $\zray \in \liets$ be such that the
distribution $\ft h$ is smooth on the ray $t \zray$ for  $t \in (0, \delta)$
for some $\d > 0 $,
and suppose $(\ft h)(t \rho) $ tends to a well defined limit
as $t \to 0^+$. Then we
define
\beq \label{8.6n}  \resrlch (h  \dps)
= \lim_{t \to 0^+} \reslch \Bigl ( h (\psi )
e^{i t \zray (\psi) }  \dps \Bigr ). \eeq
\end{definition}
Under these hypotheses, $\resrlch ( h \dps) $ is independent of $\testf$
(by Proposition \ref{p8.7n}), but it may depend on the
ray $\{ t \zray: t \in \RR^+ \}. $ However
we have  by Proposition \ref{p8.7n}
\begin{prop} \label{p8.8n} Suppose $\ft h $ is smooth at $0$. Then
the quantity $\resrlch ( h \dps)$ satisfies
$\resrlch (h \dps) = \reslch ( h \dps) $.
\end{prop}

\noindent{\em Remark:} In the light of Propositions \ref{p8.6n},
\ref{p8.7n} and \ref{p8.8n}, it makes sense to write
${\rm Res} (\Omega) $ for $\resrlch (\Omega) $ when
$\Omega = \nusym^2 h \dps $ for  a
$K$-invariant function $h: \liek \to \CC$ for which $\fk h$ is compactly
supported and $\ft (\nusym^2 h)$ is smooth at $0$. In the proof of
Theorem \ref{t8.1} which we are about to give, we shall check
the validity
of these hypotheses for the form $\Omega$ which appears in
the statement of the Theorem.

We now give the proof of Theorem \ref{t8.1}. We shall first show

\begin{prop} \label{p7.6}

{\noindent (a)} The distribution
$\fk(\pisk \eta e^{i \bom} ) (\zp) $  defined for
$\zp \in \liek$ is represented by a smooth
function for $\zp$ in a sufficiently small neighbourhood of $0$.

{\noindent (b) }  We have
\beq  \label{7.001} \eta_0   e^{i \om_0}
[\xred]  = \frac{1}{(2 \pi)^{s/2}i^s
\vol(K) } F_K (\pis \eta e^{\iins \bom}) (0),  \eeq

\beq \label{7.002} \phantom{ \eta_0   e^{i \om}
[\xred]}   = \frac{(2 \pi)^{l/2} }{(2 \pi)^{s}
|W| \vol(T)i^s } F_T (\nusym^2 \pis \eta e^{\iins \bom}) (0). \eeq
\end{prop}

\noindent{\em Proof of (a):}
 To evaluate $\fk(\pis \eta e^{i \bom} ) (0)$, we introduce
a test function $\testf: \lieks \to \RR^+$ which is smooth and
of compact support, and for which $(\fk \testf) (0) = 1/(2 \pi)^{s/2}. $
We define $\testfe(z) = \e^{-s} \testf(z/\e) $; as $\epsilon \to 0 $,
the functions $\testfe $ tend to the Dirac delta distribution on $\lieks$
(in the space $\cald'$ of distributions on $\lieks$).
Then we have
\beq (\fk \testfe) (\phi) = (\fk \testf) (\e \phi). \eeq
Now to evaluate $\fk (\pis \eta e^{i \bom} )  (\zp)$, we integrate it against
the sequence of test functions $\testfe$:
\beq \fk ( \pis \eta e^{i \bom} ) (\zp) = \lim_{\e \to 0^+}
i^s  \je (\zp) \eeq
where
\beq i^s \je (\zp)= \tintk [d z] \fk (\pis \eta
e^{i \bom} ) (z) \chi_\e (z - \zp)  \eeq
$$ \bla = \intk [d \phi] \int_{x \in X} \eta(\phi) e^{i \om}
e^{i (\mu(x) - \zp)  (\phi) }  \htestfe(\phi) $$
(by Parseval's Theorem).
Now because $\testfe$ is smooth and of compact support, the Paley-Wiener
Theorem (Theorem 7.3.1 of \cite{hor}) implies that $\htestfe$ is rapidly
decreasing. So we may use Fubini's theorem to interchange the order
of integration and get
\beq \label{7.t1} i^s \je(\zp)  = \int_{x \in X } e^{i \om} \int_{ \phi \in
\liek}
[d \phi] \eta(\phi) \htestfe(\phi) e^{i (\mu(x) - \zp) (\phi)} \eeq
\beq \label{7.t2} \bla = (2 \pi)^{s/2}
\int_{x \in X} e^{i \om} \eta( - i \partial/\partial z) \testfe(z)|_{z =
\mu(x) - \zp}. \eeq
As $\e \to 0$, $\testfe$ is supported on an arbitrarily
small neighbourhood of $\mu^{-1}(0)$. Thus, because of
Proposition  \ref{p5.2},
the integral may be replaced by an integral over $P \times \lieks$:
\beq i^s \je(z') = (2 \pi)^{s/2} \int_{(p,z) \in P \times \lieks}
e^{i \om} \eta( - i \frac{ \partial }{\partial z}) \testfe(z - \zp)
\eeq
Now by the same argument as given in the proof of Theorem \ref{t7.1}(a),
there is $\eta_0 \in \Omega^*(P/K) $ such that
\beq \eta - \pi^* \eta_0 = D \g \onebl \mbox{for some $\g$, } \eeq
where $D$ is the equivariant cohomology differential on $P \times \lieks$
and $\pi: P \times \lieks \to P \to P/K$ is the projection map.
By the argument given after (\ref{7.d3}), we have
\beq i^s \je(\zp)  - (2 \pi)^{s/2} \int_{(p,z) \in P \times \lieks}
e^{i \om} (\pi^* \eta_0) \testfe(z - \zp)  \eeq
$$ \bla = \int_{(p,z) \in P \times \lieks} e^{i \om}
\int_{ \phi \in \liek} [d \phi] D \g(\phi) \htestfe(\phi)
e^{i (z - \zp)( \phi) } $$
$$ = \int_{ \phi \in \liek} [d \phi] \htestfe( \phi)
\int_{ (p,z) \in P \times \lieks } e^{i \om} D \g (\phi)
e^{i (z - \zp) (\phi) } \eqdef \deleps. $$
But
\beq e^{i \om} e^{i (z-\zp)(\phi) }
D \g = D (e^{i \om} e^{i (z - \zp) (\phi) }
\g(\phi) ). \eeq
Hence
\beq \deleps = \intk [d \phi] \htestfe (\phi)
\int_{(p,z) \in  P \times \lieks}
d ( e^{i \om} e^{i (z - \zp) (\phi) } \g(\phi) ) \eeq
(since for  any differential form $\Psi$, the term
$\iota_{\tilde{\phi}}  \Psi$ in $D \Psi$
cannot contain differential forms of top degree in $x$).
Replacing $P \times \lieks$ by
$P \times B_R(\lieks)$ where $B_R (\lieks)$ is a large ball in $\lieks$
with boundary $S_R (\lieks)$, we have that
\beq \deleps = \lim_{R \to \infty}
\intk [d \phi] \htestfe( \phi)
\int_{(p,z) \in P \times B_R (\lieks) } d \Bigl (
e^{i \om} e^{i (z - \zp) (\phi) }
\g( \phi)    \Bigr ) . \eeq
We then interchange the order of integration to get
\beq \deleps =  \lim_{R \to \infty}
\int_{(p,z) \in P \times B_R (\lieks) } d \left (
e^{i \om} \intk [d \phi] \htestfe ( \phi) e^{i (z - \zp) (\phi) }
\g( \phi) \right ) \eeq
$$ \bla =  \lim_{R \to \infty}
\int_{(p,z) \in P \times B_R(\lieks)  } d
\Bigl (e^{i \om}
\g( \frac{\partial}{\partial z}) \testfe (z - \zp) \Bigr  ) . $$
$$ \bla =  \lim_{R \to \infty}
\int_{(p,z) \in P \times S_R(\lieks)  }
\Bigl (e^{i \om}
\g( \frac{\partial}{\partial z}) \testfe (z - \zp) \Bigr  )
\onebl \mbox{by Stokes' Theorem}. $$
This limit
equals $0$ for sufficiently
small $\zp$ since $\testfe$ is
compactly supported on $\lieks$.  Hence $\deleps = 0 $.

Finally we obtain using the expression for $\omega$
given in Proposition \ref{p5.2}
\beq i^s \je(\zp)  = (2 \pi)^{s/2} \int_{(p,z) \in P \times \lieks} e^{i \om}
(\pi^* \eta_0) \testfe (z - \zp), \eeq
\beq \bla = ( 2 \pi)^{s/2} \int_{ (p,z) \in P \times \lieks}
e^{i \pi^* \om_0} (\pi^* \eta_0) e^{i \theta(dz)} e^{i d \theta(z) }
\testfe (z - \zp). \eeq
Thus we have as in Section 5
\beq \label{8.26}
 \lim_{\e \to 0 } \je (\zp) = ( 2 \pi)^{s/2} \int_{(p,z) \in
P \times \lieks} e^{i \pi^* \om_0} (\pi^* \eta_0)
\Om [dz] e^{i F_\theta (z)} \delta(z - \zp) \eeq
$$ = ( 2\pi)^{s/2} \int_P e^{i \pi^* \om_0}
 (\pi^* \eta_0)
e^{i F_\theta(\zp) }
\Om $$
where $\Omega$ is the differential form introduced after (\ref{5.4}).
This shows in particular that $(\fk  \rbare)(\zp)$
 (where $\rbare (\phi) = \pis (\eta e^{i \bom})( \phi) $)
is a
polynomial in $\zp$ for small $\zp$, and hence is smooth  in
$\zp$ for $\zp$ sufficiently close to $0$. This completes the
proof of (a).

When $\zp = 0 $,  equation (\ref{8.26})
becomes
$$\lim_{\e \to 0}
\je (0) =
i^{-s} F_K (\pis \eta e^{\iins \bom})(0) = {(2 \pi)^{s/2}   \vol(K) }
\eta_0 e^{i \om_0} [\xred], $$
which proves (\ref{7.001}).
Using Lemma \ref{l3.1} again, we have for any
$K$-invariant function $f$ on $\liek$ that
\beq \label{7.003}
\frac{(2 \pi)^{s/2} }{ \vol(K)  }
(F_K f) (0) =
\frac{( 2 \pi)^{l/2} }{|W| \vol(T) }
F_T(f \nusym^2) (0). \eeq
Combining this with (\ref{7.001}) we obtain
$$
\eta_0 e^{i \om_0}[\xred]
   = \frac{(2 \pi)^{l/2} }{(2 \pi)^{s}
|W| \vol(T) i^s  } F_T (\nusym^2 \pis \eta e^{\iins \bom}) (0), $$
which is (\ref{7.002}). $\square$

Theorem \ref{t8.1} follows from Proposition \ref{p7.6} and
Proposition \ref{p8.7n}
by applying
Theorem \ref{t2.1}  to decompose $r^{{\eta}}  =
\pis (\eta e^{i \bom} ) $ as a sum
of meromorphic functions $\rfe$ on $\liet \otimes \CC $ corresponding
to the components $F$ of the fixed point set of $T$.  We now
complete the proof of this theorem.

\noindent{\em Proof of Theorem \ref{t8.1}:} As in (\ref{7.6n}),
the abelian localization formula yields
$ r^{{\eta} } (\psi) = \sum_{F \in \calf}  \rfe (\psi), $
where
$$ \rfe(\psi) = e^{i \mu_T(F) (\psi) } \int_F
\frac{i_F^* \eta(\psi) e^{i \om} } {e_F(\psi) }. $$
Now the distribution $(\fk r^{{\eta} } ) (\phi) $ is represented
by a smooth function near $0$ (Proposition \ref{p7.6} (a));
also, the distribution $\ft(\nusym^2 r^{{\eta} } ) $
$ = D_{\nusym} \ft (\nusym \rbare) $ is represented
by a smooth function near $0$ since $\ft (\nusym r^{{\eta} } )  =
\nusym \fk r^{{\eta} } $ (Lemma
\ref{l3.3}) and $\fk r^{{\eta} } $ is smooth near $0$ (Proposition
\ref{p7.6}(a)).
We choose $ \zray \in \liets$ so that the distribution
$\ft \rfe$ is smooth along the ray $t \zray, t \in (0, \delta) $, for all
$F$ and sufficiently small $\delta > 0 $, and that this distribution
tends to a well defined limit as $t \to 0^+$:
this is possible because the $\rfe (\psi)$ are sums of
terms of the form $e^{i \mu_T(F) (\psi) } / \prod_j \bfj(\psi)^{n_j} $
so their Fourier transforms $\ft \rfe$ are piecewise polynomial
functions of the form $H_\barb (y) $ (see Proposition \ref{p3.5}).
These functions are smooth on the set $U_\barb$ consisting of all points
$y$ where $y$ is not in the cone spanned by any subset of the
$\bfjw$ containing less than $l$ elements. Thus
$$\lim_{t \to 0^+} (2 \pi)^{-l/2}  i^{-l}   \ft \rfe(t \zray) =
 \resrlch (\rfe \dps) $$ by
Definition \ref{d8.5n} and Proposition \ref{p8.7n}. It follows that
\beq \label{8.7n}
(2 \pi)^{-l/2} i^{-l}  \ft (\nusym^2 \rbare) (0) =
\resrlch (\nusym^2 \rbare \dps ) \eeq
\beq  \label{8.8n} \bla  = \sum_{F \in \calf} \resrlch (\nusym^2 \rfe).  \eeq
The residue in (\ref{8.7n}) is independent of $\testf$, $\lasub$ and $\zray$
by Propositions \ref{p8.6n}, \ref{p8.7n} and
\ref{p8.8n}. The residues in (\ref{8.8n}) are independent of $\testf$ by
Proposition \ref{p8.7n}, but they may depend on $\zray$
and $\lasub$.

To conclude the proof of Theorem \ref{t8.1} we note that (\ref{7.002})
gives
\beq
  \eta_0   e^{i \om}
[\xred]   = \frac{(2 \pi)^{l/2} }{(2 \pi)^{s}
|W| \vol(T)i^s } F_T (\nusym^2 \pis \eta e^{\iins \bom}) (0),  \eeq
$$   = \frac{i^l  }{(2 \pi)^{s-l}
|W| \vol(T)i^s }  \reso(\rbare(\psi)    \nusym^2(\psi)   [ d\psi]) \onebl
\mbox{(by Proposition \ref{p8.7n})} $$
$$  = \frac{-1)^\npl }{ (2 \pi)^{s-l}  |W|  \vol(T) }   \reso
\Biggl (  \sum_{F \in \calf} \rfe(\psi)  \nusym^2 (\psi) [d \psi]
\Biggr ) $$
as claimed. $\square$

\noindent{\em Proof of Corollary 8.2:} In a normalization
where $\vol(T) = 1$, the factor $\nusym(\psi) $ becomes $2 \pi \psi$.
This gives
$$ \frac{1}{(2 \pi)^{s-l}}\treso
\Biggl (     \: \nusym^2 (\psi)
 \sum_{F \in \calf} \rfe(\psi) \dps \Biggr )
= \reslch  (  \psi^2
 \sum_{F \in \calf} \rfe(\psi) [d \psi] ). $$
Each term $\rfe (\psi) $ is a sum of terms of the form
$\tau_\a (\psi) = c_{\a} \psi^{-n_\a}  e^{i \mu(F) \psi}  $ for some constants
$c_\a$ and integers $n_\a$.
By (\ref{8.1n}), the residue is given by
\beq \label{8.001} \reslch (h \dps) =
\lim_{ \e \to 0^+ }
\frac{1} {2 \pi i} \int_{ \psi \in \RR - i \xi}
\htestf (\e \psi) h(\psi) \dps, \eeq
where we choose $\xi$ to be in the cone $\lasub = \RR^+$.
Proposition \ref{cty} (i), (ii) shows that
the function $\htestf: \RR \to \CC$ extends to an entire
function on $\CC$.

We may now  decompose the integral
 in (\ref{8.001}) into terms corresponding to the $\tau_\a$. If
$n_\a > 0 $, we
complete each such integral over $\RR - i \xi$ to a
contour integral by adding a semicircular curve  at infinity, which
is in the upper half plane if $\mu(F) > 0 $ and in the lower half
plane if $\mu(F) < 0 $. This choice of contour is  made so that
the function
$\tau_\a (\psi) $  is bounded on the added
contours, so the added semicircular curves  do not contribute to
the integral. Since only the contours corresponding to values
of $F$ for which $\mu(F) > 0 $  enclose the pole at $0$, application
 of Cauchy's residue formula  now gives the result.
A similar argument establishes that the terms  $\tau_\a $ for which
$n_\a \le  0 $ contribute $0$ to the sum.\footnote{The
formula obtained by choosing $\lasub = \RR^-$ is in fact
equivalent to the formula we have
obtained using the choice $\lasub = \RR^+$.
 This can be seen directly from the Weyl invariance
of the function $\rbare$, where the action of the Weyl group takes
$\psi$ to $- \psi$ and so converts terms with $\mu(F) > 0  $ to
terms with $\mu(F) < 0 $.}
$\square$

\noindent{\em Remarks:}

 (a) The quantity $\resrlch  (\nusym^2 \rfe \dps ) $ depends on
the cone $\lasub$ for each $F$;
however, it follows from Proposition
\ref{p8.6n} that the sum $\sum_{F \in \calf}
\resrlch  ( \nusym^2 \rfe \dps ) $ is independent of $\lasub$. }

\noindent (b) Let $\calf_\lasub$ be the set of those $F \in \calf$ for which
$\mu_T(F)$ lies in the cone $\dcf =
\{ \sum_j s_j \bfjw: s_j \ge 0 \} $ (defined in  Section 4) spanned
by the $\bfjw$. Then by
Proposition \ref{p8.4} (iii), $\resrlch ( \nusym^2 \rfe ) = 0 $ if
 $F \notin \calf_\lasub$,
so in fact $\resrlch (\nusym^2  \rbare \dps)  =
\resrlch \sum_{F \in \calf_\lasub} ( \nusym^2 \rfe  \dps ) $.

\noindent (c)  Finally, it follows from Proposition \ref{p8.4} that
if we replace the symplectic form
$\omega$ and the moment map $\mu$ by $\delta \omega$ and $\delta \mu$
(where $\delta > 0 $) and then let $\delta $ tend to $0$, we obtain
an expression where $\mu$ and $\omega $ appear only in determining the
set $\calf_\lasub$ indexing terms which yield a nonzero contribution.

We now restrict ourselves to the special case when $\Omega$ is of the
form $$\Omega_\l(\psi)  = e^{i \l(\psi)}
[d \psi]/\prod_{j = 1}^N \beta_j (\psi).$$
If $\b_j \in \lasub $ then the distribution
$\resrlch (\Omega_\l) $ is just (up to multiplication by a
constant) the piecewise polynomial function $H_\barb(\l)$ from Proposition
\ref{p3.5}.
We shall now give  a Proposition
 which gives a list of properties satisfied
by the residues $\reso (\Oma{\l})$: these properties in fact characterize the
residues uniquely and enable one to compute them.

\nc{\pee}{P}
\begin{prop} \label{p8.4} Let $\xi \in \liet$ and
suppose
 $\b_1, \dots, \b_N \in \liets$ are
all in the dual cone of a cone $\lasub \in
\liet$.
Denote by $\pee: \liet \to \RR $ the
function $\pee(\psi) = \prod_j \b_j(\psi)$.
Suppose  $\l \in U_\barb \subset \liets$ (see
Proposition \ref{p3.5}), and
define $\Oma{\l}(\psi) =
e^{i \l(\psi)}[d \psi]/\pee(\psi) $.

Then we have
\begin{description}
\item[(i)]
 $\reso (\psi_k \psi^J \Oma{\l})
= (- i \partial/\partial \lambda_k) \reso (\psi^J \Oma{\l}).
$
\item[(ii)] $(2 \pi i)^l \reso (\Oma{\l})
=  i^N H_\barb (\lambda)$, where $H_\barb$ is
the distribution given in Proposition \ref{p3.5}. (Recall we have assumed
that
$\b_j $ is in the dual cone of
$ \lasub $ for all $j$.)
\item[(iii)] $\reso (\Oma{\l}) = 0 $ unless $\lambda$ is in the cone
$C_\barb$ spanned
by the $\beta_j$.
\item[(iv)] $$\lim_{s \to 0^+} \reso \Bigl (
\Oma{s \l} \psi^J
\Bigr ) = 0$$
unless $ N  - |J| =  l$.
\item[(v)]
$$\lim_{s \to 0^+} \reso \Bigl (
\Oma{s \l} \psi^J
\Bigr ) = 0 $$
if the monomials $\beta_j$ do not span $\liets$.
\item[(vi)] If
$\b_1, \dots, \b_l$ span $\liets$ and
$\lambda = \sum_j \lambda^j \b_j$
with all $\lambda^j > 0$,
 then
$$ \lim_{s \to 0^+} \reso
\Bigl ( \frac{ e^{i s \lambda(\psi) } [d \psi] }{\beta_1(\psi) \dots
\beta_l(\psi)
}   \Bigr )
= \frac{1}{\det \barb } ,  $$
where $\det \barb $ is the determinant of the $l $ by $l$
matrix whose columns are the coordinates of $\beta_1, \dots,
\beta_l$ written in terms of any orthonormal
basis of $\liet$.
\item[(vii)] $$ \reso  \Bigl (
\frac{e^{i \lambda(\psi) }  \psi^J   }
{\pee(\psi)}
[d \psi]
\Bigr  ) = \sum_{m \ge 0} \lim_{s \to 0^+}
\reso \Bigl (  \frac{ (i \lambda(\psi))^m e^{i s \lambda(\psi) } \psi^J  }
{m! \pee(\psi)}
 [d \psi]\Bigr ). $$

\end{description}
\end{prop}

\noindent{\em Remark:} The limits in Proposition
\ref{p8.4} are not part of the
definition of the residue map: rather these limits
are described in order to specify a procedure
for computing the piecewise polynomial function
$H_\barb (\l) = (2 \pi i)^l i^{-N}
\reso (\Om_\l). $ (See the example below.)
Proposition \ref{p8.4} (ii)
identifies
$H_\barb (\l)$ with  an integral over $\liet$,
which may be completed to an appropriate
contour integral: the choice of contour is
determined by the value of $\l$,
and requires $\l$ to be a nonzero point in
$U_\barb$. Proposition \ref{p8.4}
(vii) says that one may compute $H_\barb$ by
expanding the numerator $e^{i \l(\psi)}$ in
a power series, but only provided one keeps
a factor $e^{is \l(\psi)}$ in the
integrand (for small $s > 0$) in order to
specify the contour. The limits in
Proposition \ref{p8.4} (iv)-(vii) exist because
according to Proposition \ref{p8.4}
(i) and (ii), they specify limits
of derivatives of polynomials on subdomains of
$U_\barb$, as one approaches the point $0$
in the boundary of $U_\barb$ along the fixed
direction $s \l$ as $s \to 0$ in $\RR^+$.

\noindent{\em Proof of  (i):} This
follows directly from Definition \ref{d8.5n}.

\noindent{\em Proof of (ii):}
This follows because $(2 \pi i)^l i^{- N} \reso(\Oma{\l})$
is the fundamental solution $E(\l)$ of the differential
equation $P(\partial/\partial \l)E(\l) = \delta_0$
with support in a  half space containing the $\b_j$.
(See
\cite{abg} Theorem 4.1 or  \cite{hor2} Theorem 12.5.1.)
But this fundamental solution is given
by $H_\barb$ (see Proposition  \ref{p3.5}(c)).

\noindent{\em Proof of (iii):}
This is an immediate consequence of (i) and (ii), in view
of Proposition  \ref{p3.5}(a).\footnote{Alternatively there is the
following direct argument, which was pointed out to us
by J.J. Duistermaat. We recall that
$$
 \reso (\Omega_\l) = \frac{1}{(2 \pi i)^l} \intt
\frac{ e^{i \l(\psi - i \xi) }}{P(\psi - i \xi) }  \, [d \psi], $$
which is defined and independent of $\xi$ for  $\xi \in (C_\barb)^* $
(see the discussion after
Definition \ref{d8.5n}). Hence we may replace $\xi $ by $t \xi $
for any $t \in \RR^+$:
\beq \label{resrec}
 \reso (\Omega_\l) = \frac{1}{(2 \pi i)^l} \intt
\frac{ e^{i \l(\psi - i t\xi) }}{P(\psi - i t \xi) }  \, [d \psi].  \eeq
Taking the limit as $t \to \infty$, we
see that $\reso (\Omega_\l) = 0 $ if $ \l(\xi) < 0 $,
because of the factor $e^{  t \l(\xi) } $ that appears
in the numerator of (\ref{resrec}). Since this holds
for all $\xi \in (C_\barb)^*$,
$\reso (\Omega_\l)$ is only nonzero when $\l(\xi) \ge 0 $
for all $\xi \in (C_\barb)^*$, in other words when
$\l \in C_\barb$.}

\noindent{\em Proof of (iv):} By (i),
$$
\limpl \reso \Bigl ( \Oma{s \l_0} \psi^J
 \Bigr ) = \limpl (- i \partial/\partial \l)^J
H_\barb (\l)|_{\l = s\l_0}, $$
but $H_\barb$ is a homogeneous piecewise polynomial function
of degree $N-l$, so the conclusion holds for $|J| > N - l.$
If $|J| < N - l$, we find that $(\partial/\partial \l)^J
H_\barb(\l)$ is homogeneous of order $N - l - |J|$ (on any
open  subset of
$\liets$ where $H_\barb$ is smooth). Hence it
is of order $s^{N - l - |J|} $ at $\l = s \l_0$ as $s \to 0^+$,
and the conclusion also holds in this case.

\noindent{\em Proof of (v):} By (ii), we know that
$\reso \Bigl ( \Oma{\l}
 \Bigr )
 = 0$ for $\l$ in a neighbourhood of $s \lambda_0$ (since $s \lambda_0$
is not in the support of $H_\barb$). Applying
(i), $\reso  \Bigl (
\Oma{\l} \psi^J
\Bigr ) $
must also be zero.

\noindent{\em Proof of (vi):} If $\b = \{\b_1, \dots,
\b_l \} $ span $\liets$, and $\lambda = \sum_j \lambda^j \b_j$
with all $\lambda^j > 0$, we have
$$\reso \Bigl ( \Oma{\l} \Bigr )  = \frac{1}{(2 \pi i)^l}
\intt \frac{ [d \psi] e^{i \sum_j \lambda^j (\psi_j - i  \vt_j) } }
{\prod_{k = 1}^l (\psi_k - i \vt_k ) } , $$
where the $\vt_k = \b_k (\vt) > 0$.
But $ [d \psi]  = d \psi_1 \dots d \psi_l/(\det \barb)$,
where $\psi_j = \b_j(\psi)$.
Since the integrals over $\psi_1, \dots, \psi_l$ may be
completed to integrals over semicircular contours
$C_+(R) $  in the
upper half plane, for each of which the contour integral
may be evaluated by the Residue Theorem to give the
contribution $ 2 \pi i$, we obtain the result.

\noindent{\em Proof of (vii):} We have
$$ \int_{\psi + i \vt \in \RR^l} \frac{  e^{i \l(\psi)} }{P(i\psi)} [d \psi]
= (2 \pi i)^l i^{-N} \reso(\Oma{\l}) =
H_\barb (\l). $$
Also, by (i),
$$\sum_{m_1, \dots,
m_l \ge 0}  \frac{(2 \pi i)^l i^{-N} }{m_1!\dots m_l! }
\limpl
\reso \Bigl (
(i \l^1\psi_1)^{m_1}
\dots (i \l^l\psi_l)^{m_l} \Oma{s \l_0} \Bigr )  $$
$$ = \sum_{m_1, \dots,
m_l \ge 0}  \frac{(2 \pi i)^l i^{-N}
(\l^1)^{m_1}  \dots
(\l^l)^{m_l}}
{m_1!\dots m_l! }
\limpl ( \partial/\partial \l^1)^{m_1} \dots
( \partial/\partial \l^l)^{m_l}  \reso (\Oma{\l} )
|_{\l = s \l_0} $$
$$ = \sum_{m_1, \dots,
m_l \ge 0}  \frac{
(\l^1)^{m_1}  \dots
(\l^l)^{m_l}}
{m_1!\dots m_l! }
\limpl ( \partial/\partial \l^1)^{m_1} \dots
( \partial/\partial \l^l)^{m_l}  H_\barb(\l)
|_{\l = s \l_0} $$
which is equal to $H_\barb(\l)$ since $H_\barb$ is
a polynomial on certain conical subregions
of $\liets$, and there is such
a subregion containing the ray $\l = s \l_0$.$\square$

\noindent{\underline{ Example.}}
In the following  simple example, the explicit formula
 for $H_\barb$ follows immediately from the definition of $H_\barb$ in
Proposition \ref{p3.5}(a). The example is included to show how
this result may alternatively be derived by successive contour
integrations.

Suppose $l = 2$ and
$N = 3$, and $\b_1(\psi) = \psi_1$, $\b_2(\psi) = \psi_2$,
$\b_3(\psi) = \psi_1 + \psi_2$. We compute
\beq \label{10.6}\calr = \reso \Bigl ( \frac{e^{i \l(\psi)} }
{\psi_1 \psi_2 (\psi_1 + \psi_2)} \Bigr ) \eeq
where $\l(\psi) = \l^1 \psi_1 + \l^2 \psi_2. $
We assume $\l^1, \l^2 > 0$, and $\vt_1, \vt_2 > 0 $.
The quantity (\ref{10.6}) is given by
\beq \calr  =
\frac{1}{(2 \pi i)^2}
\int_{\psi_2 + i \vt_2 \in \RR} d \psi_2 \frac{e^{i \l^2 \psi_2}}{\psi_2}
\int_{\psi_1+ i \vt_1\in \RR}
 \frac{e^{i \l^1 \psi_1}}{(\psi_1 + \psi_2)\psi_1} . \eeq
(This integral in fact gives the Duistermaat-Heckman
polynomial $F_T (\pist e^{\iins \bom}) (\l)$ near $\mu_T(F) $ where
$F$ is
  a fixed point of the action of
$T$ on $X$,  when $X$ is a
coadjoint orbit of $SU(3)$ and $T\cong (S^1)^2$ is the
maximal torus : see \cite{JGP}.)
We compute this by first integrating over $\psi_1$: since $\l^1 > 0 $,
the integral may be completed to a contour integral over a semicircular
contour in the upper half plane. We obtain contributions
from the two residues $\psi_1 = 0 $ and $\psi_1 = - \psi_2$.
Hence we have
\beq \calr =
\frac{1}{(2 \pi i)}
\int_{\psi_2 + i \vt_2 \in \RR} d \psi_2 \frac{e^{i \l^2 \psi_2}}{\psi_2^2}
-\frac{1}{(2 \pi i)}
\int_{\psi_2 + i \vt_2 \in \RR} d \psi_2 \frac{e^{i (\l^2-
\l^1) \psi_2}}{\psi_2^2} . \eeq
Since $\l^2 > 0 $, the first of these integrals may
be completed to a contour integral over a semicircular
contour in the upper half plane, and the residue
at $0$ yields the value $i \l^2$. If $\l^2 - \l^1 > 0 $, the second
integral likewise yields $- i (\l^2 - \l^1)$.
However if $\l^2 - \l^1 < 0 $ the second integral
is instead equal to a contour integral over a semicircular
contour in the {\em lower} half plane, which
does not enclose the pole at $0$, and hence the second integral gives $0$.
Thus we have
\beq \label{9.009}(2 \pi i)^2 \reso(\Oma{\l}) = \cases{i \l^2, &
$\l^1 > \l^2$ \cr
i \l^1, &
$\l^2 > \l^1$. \cr } \eeq

According to (iv) and (vii), the quantity $\calr$ is
also given by
\beq \calr = \limpl \int_{\psi+ i \vt\in \RR^2}
\frac{ e^{i s (\l^1 \psi_1 + \l^2 \psi_2)} (i \l^1 \psi_1
+ i \l^2 \psi_2)[d \psi]}{\psi_1 \psi_2 (\psi_1 + \psi_2)} \eeq
$$ =
  \limpl (i \l^1) \int_{\psi+ i \vt\in \RR^2}
\frac{ e^{i s [\l^1 (\psi_1 + \psi_2) +  (\l^2 - \l^1)\psi_2]}
[d \psi] }{(\psi_1 + \psi_2) \psi_2}
+
  \limpl (i \l^2) \int_{\psi+ i \vt\in \RR^2}
\frac{ e^{i s [(\l^1 - \l^2)\psi_1 +  \l^2 (\psi_1 + \psi_2)]}
[d \psi] }{\psi_1 (\psi_1 + \psi_2) }
. $$
This clearly gives the result  (\ref{9.009}).

\renorm
\section{Examples}

In this section we shall show
in the case $K = SU(2)$  how Corollary  \ref{c8.2} may
be used to prove relations in the cohomology
ring $H^*(\xred)$ for two specific $X$.
These $X$ are the examples treated at the end of
Section 6 of \cite{Ki2}. There, all the relations
in the cohomology ring are determined, which is equivalent
to exhibiting
all the   vanishing intersection
pairings. We shall show how the results
of the present paper may be used to show these are indeed vanishing
  intersection pairings, although we shall not
be able to rederive the result that there are no others.

\noindent\underline{\em Example 1: $X = ({\Bbb P}_1)^N, $  $N$ odd.}
Consider the action of $K = SU(2)$ on the space
$X = (\PP_1)^N$ of ordered $N$-tuples
of points on the complex projective line
$\PP_1$, defined by the $N$th tensor
power of the standard representation of $K$ on
$\CC^2$.
Equivalently when $\PP_1$ is identified with
the unit sphere $S^2$ in $\RR^3$ then $K$ acts
on $X = (S^2)^N$ by rotations of the sphere. When
the dual of the Lie algebra of $K$ is identified
suitably with $\RR^3$ then the moment map
$\mu$ is given (up to a constant scalar
factor depending on the conventions used)
by $$\mu(x_1, \dots, x_N)  = x_1 + \dots + x_N$$
for $x_1 , \dots, x_N \in S^2$. We assume
that $0$ is a regular value for $\mu$; this happens
if and only if there is no $N$-tuple
in $\mu^{-1}(0)$ containing a pair of antipodal
points in $S^2$ each with multiplicity $N/2$,
and so $0$ is a regular value if and only
if $N$ is odd.

In order to apply Corollary \ref{c8.2} we note that
the fixed points of the action of the standard
maximal torus $T$ of $K$ are the $N$-tuples
$(x_1, \dots, x_N)$ of points in
$\PP_1$ such that each
$x_j$ is either $0$ or $\infty$. We shall index these
by sequences $n = (n_1, \dots, n_N)$  where
$n_j = + 1$ if $x_j = 0$ and $n_j = - 1$
if $x_j = \infty$.
Denote by $e_n$ the fixed point indexed
by $n$. Then
$$ \mu_T(e_n) = \sum_{j = 1}^N n_j $$
and the weights of the action of $T$
at $e_n$ are just $\{ n_1, \dots, n_N\}$.
Hence the sign of the product of weights
at $e_n$ is $\prod_j n_j$ and its
absolute value is $1$.

The cohomology ring $H^*(X)$ has $N$
generators $\xi_1, \dots, \xi_N$ say, of degree two,
satisfying $\xi_j^2  = 0$
for $1 \le j \le N$. The equivariant cohomology
ring $H^*_T(X)$ with respect to the torus $T$
has generators
$\xi_1, \dots, \xi_N$ and $\a$ of
degree two subject to the relations
$$(\xi_j)^2 = \a^2$$
for $1 \le j \le N$. The Weyl group action sends $\a$
to $-\a$ so $\hk(X)$ has generators
$\xi_1, \dots, \xi_n, \a^2$ subject to the
same relations.

According to the last example of section 6 of
\cite{Ki2}, the kernel of the map $\hk(X) \to H^*(\xred)$
is spanned by elements of the form
\beq \label{8.beta} (1/\a) \Bigl ( q(
\xi_1, \dots, \xi_N, \a) \prod_{i \in Q} (\xi_i + \a)
- q(
\xi_1, \dots, \xi_N, -\a) \prod_{i \in Q} (\xi_i - \a) \Bigr )  \eeq
for some $Q \subset \{ 1, \dots, N\} $  containing
at least $(N+1)/2$ elements
and some polynomial $q$ in $N + 1$ variables
with complex coefficients.
We can  use Corollary \ref{c8.2}
to give an alternative proof that the evaluation against the
fundamental class $[\xred]$ of the image in
$H^*(\xred)$ of any element
of this form of degree $N-3$ is zero. This amounts
to showing that
\beq {\rm Res}_0 \frac{1}{\psi^{N-1} }
\sum_{\stackrel{\indd_j = \pm 1,}{ \sum_j \indd_j > 0} } \Bigl ( \prod_j
\indd_j  \Bigr )
\Bigl \{ q(\psi \indd_1, \dots, \psi \indd_N, \psi) \prod_{i \in Q}
\psi(\indd_i + 1)
- q(\psi \indd_1, \dots, \psi \indd_N, -\psi) \prod_{i \in Q}
\psi(\indd_i - 1) \Bigr \} \eeq
is zero when $q$ is homogeneous of degree $N - 2 - |Q|$.
In other words it amounts to showing
that $\dell = 0$ for any $\dell$ of the form
\beq \dell =
\sum_{\stackrel{\indd_j = \pm 1,}{ \sum_j \indd_j > 0}}
\Bigl ( \prod_j \indd_j  \Bigr )
\Bigl \{
q( \indd_1, \dots,  \indd_N, 1) \prod_{i \in Q}
(\indd_i + 1)
- q( \indd_1, \dots,  \indd_N, -1) \prod_{i \in Q}
(\indd_i - 1) \Bigr \} \eeq
where $q$ is homogeneous of degree $N - 2 - |Q|$.
Let us assume without loss of generality that
$q (\xi_1, \dots, \xi_N, \a) = \prod_i \xi_i^{r_i} \a^p$
where $p + \sum_i r_i + |Q| = N-2$.
Thus $p - 1 = |Q| + \sum_i r_i $ (mod $2$) since
$N$ is odd. Hence we have that
\beq \dell =
 \sum_{\stackrel{\indd_j = \pm 1,}  {\sum_j \indd_j > 0}
}
\prod_j  (\indd_j^{r_j + 1} ) \Bigl \{
\prod_{i \in Q} (\indd_i + 1) +
(-1)^{|Q|} \prod_k (-1)^{r_k} \prod_{i \in Q} (\indd_i - 1)  \:
\Bigr \}  \eeq
\beq \onebl = \sum_{\stackrel{\indd_j = \pm 1,}  {\sum_j \indd_j > 0}
}
\prod_j  (\indd_j^{r_j + 1} )
\prod_{k \in Q} (\indd_k + 1)
-
 \sum_{\stackrel{\Indd_j = \pm 1,}{ \sum_j \Indd_j < 0} }
\prod_j  (\Indd_j^{r_j + 1} )
\prod_{k \in Q} (\Indd_k + 1)  \eeq
where we have introduced $\Indd_j = - \indd_j$.

Now the second sum vanishes, for if $\Indd_j = 1$
for all $j \in Q$ then we must have $\sum_j \Indd_j > 0$
since $|Q|> N/2$.
Hence we are reduced to proving the vanishing of
$$ \sum_{  \indd\in \G}
\prod_j \indd_j^{r_j + 1} $$
where
$$\G = \{ \indd \: | \: \sum_j \indd_j > 0, \;
\indd_j = 1 \; \mbox{for $j \in Q$} \}. $$
Hence we have to prove the vanishing
of $$\dell = \sum_{\indd_j = \pm 1, j \notin Q}
\prod_{j \in S} n_j$$
where $S$ is the set $\{ j \notin Q \, | \,
r_j = 0 \pmod{2} \}. $
The sum thus vanishes by cancellation in pairs provided $S$ is
nonempty. However if $S$ were empty then
$r_j = 1 \pmod{2}$ for all $j \notin Q$,
so $r_j \ge 1$ for all $j \notin Q$,
which is impossible since $\sum_j r_j + |Q| \le N-2$.
This proves the desired result.

\noindent\underline{\em Example 2: $X = \PP_N$, $N$ odd.}
 A closely related example
is given by the action of $K = SU(2)$ on the complex projective
space
$X = {\Bbb P}^N $  defined by the $N$th symmetric
power of the standard representation of $K$ on
$\CC^2$. Equivalently we can identify $X$ with
the space of
of {\em unordered}
$N$-tuples of points in  the complex
projective line $\PP_1$ or the sphere $S^2$, and
then $K$ acts by rotations as in example $1$.
We take the symplectic form $\om$ on $X$
to be the Fubini-Study form on $\PP_N$.
The moment map is given  by the composition
of the restriction map ${\bf u}(N+1)^* \to \lieks$
with the map
$\mu: \PP_N \to {\bf u}(N+1)^* $ defined for
$a \in {\bf u}(N+1)$ by
$$ \inpr{\mu(x) , a} = (2 \pi i |x^*|^2 )^{-1}
\bar{x^*}^t a x^*, $$
where $x^* = (x^*_0, \dots, x^*_N)$ is
any point
in $\CC^{N+1}$  lying over the point $x  \in \PP_N$.
The restriction of this moment map
to $\liets$ is
$\mu_T (x) = (2 \pi i |x^*|^2 )^{-1} )
\sum_{j = 0}^N (N - 2j)|x^*_j|^2  $.
Again we assume that $N$ is odd in order
to ensure that $0$ is a regular value of
the moment map $\mu$.

Again the fixed points of the action of $T$ are
the $N$-tuples of points in $\PP_1$ consisting
entirely of copies of $0$ and $\infty$. Equivalently
they are the points
$e_0 = [1, 0, \dots, 0]$,
$e_1 = [0, 1, \dots, 0]$,
$\dots, e_N = [0, \dots, 0, 1]$ of
$\PP_N$. The image of $e_k$ under $\mu_T$ is
$\mu_T(e_k) = N - 2 k = \mu_k$ say. Since
$\diag (t, t^{-1} ) \in T$ acts on $\PP_N$
by sending $[x_0, \dots, x_j, \dots, x_N]$
to $[t^{-N} x_0, \dots, t^{2j-N} x_j, \dots, t^N x_N]$
the weights at $e_k$ are
$$\{ 2 (j - k): \: 0 \le j \le N, \; j \ne k \}$$
The number of negative weights at $e_k$ is equal to $k$
(modulo $2$), and the absolute value of the product
of weights at $e_k$ is
$v_k = \prod_{j \ne k } |j - k| = 2^N k! (N-k)! $

The cohomology ring $H^*(\PP_N)$ is generated
by $\xi$ of degree two subject to the relation
$\xi^{N+1} = 0$.
The equivariant cohomology ring $H^*_T (\PP_N)$
is generated by $\xi$ and $\a$ of degree two
subject to the
relation $\prod_{0 \le j \le N} (\xi - (2j-N) \a) = 0, $
and the equivariant cohomology ring $\hk(\PP_N)$
is generated by $\xi$ and $\a^2$ subject to the
same relation.

According to section 6 of \cite{Ki2} the kernel
of the natural map $\hk(X) \to H^*(\xred)$
is generated as an ideal in $\hk(X) $ by $P_+(\xi, \a)$
and $P_-(\xi,
\a)/\a$ where
$$P(\xi, \a) = \prod_{k > N/2} (\xi + \mu_k \a) $$
and
$$P_\pm(\xi, \a) = P (\xi, \a) \pm P(\xi, - \a). $$
(Note that  $P_+ (\xi, \a)$ and $P_-(\xi, \a)/\a$
are actually polynomials in $\xi$ and $\a^2$.)
We would like to check that the evaluation against the
fundamental class $[\xred]$
of the image of $R_+(\xi, \a^2) P_+(\xi, \a) $
and $R_-(\xi, \a^2)  P_-(\xi, \a)/\a$ in
$H^*(\xred)$ is zero for any $R_\pm(\xi, \a^2)$
$\in \hk(X)$ of the appropriate degree.

Now we have from the abelian fixed point formula that  for any
$S(\gnx, \gna^2)$,
$$\Pi^+_*  (S(\gnx, \gna^2) ) =
\sum_{k < N/2}  (-1)^k
\frac{S(\mu_k \psi, \psi^2) }{\veee_k  }
\psi^{-N}. $$
(Here, if  $\z \in \hk(X)$, the notation $\Pi^+(\z)$
means  the portion of the abelian formula (\ref{2.1})
for $\pis(\z)$
corresponding to fixed points $F$ for which $\mu_T(F) > 0$.)
To evaluate this on the fundamental class of $\xred$
we must then find the term of degree $-1$ in
$\nusym^2(\psi) \Pi^+_* (S(\gnx, \gna^2 ) )(\psi) $, or in other words
the term of degree $N - 3$ in
$\sum_{k < N/2}(-1)^k  S(\mu_k \psi, \psi^2) /(\veee_k )
$.
Having found the term of degree $N-3$ in $\psi$, we evaluate
it at $\psi = 1$ to get the residue.
In the case when $S(\mu_k \psi, \psi^2) = R_+(\mu_k \psi, \psi^2)
P_+(\mu_k \psi, \psi) $ or
$S(\mu_k \psi, \psi^2) = R_-(\mu_k \psi, \psi^2)
P_-(\mu_k \psi, \psi)/\psi $
is homogeneous of  degree $N - 3$ in $\psi$, we need to show that
$$ \sum_{k = 0 }^{(N-1)/2}   \:  (-1)^k
\frac{1}{\veee_k } \Bigl ( \,
R_\pm (\mu_k, 1) \prod_{j \ge N/2} (\mu_k + \mu_j)
\pm  R_\pm(\mu_k, -1) \prod_{j \ge N/2} (\mu_k - \mu_j) \: \Bigr )  = 0. $$
Since $\mu_k = - \mu_{N-k} $, we have $\prod_{j \ge N/2}
 (\mu_k + \mu_j ) = 0$
for all $k < N/2$, so we just have to prove the vanishing of
$$\sum_{k \le  N/2}  \: (-1)^k
\frac{1}{ \veee_k }
R(\mu_k) \prod_{j \ge N/2} (\mu_k - \mu_j)
 $$
for every polynomial $R$ of degree at most $(N-1)/2 - 2$,
or without loss of generality the vanishing of
$$ \sum_{k \le N/2} (\mu_k)^s \frac{(-1)^k}{ \prod_{l \ne k}
|l - k| } \prod_{j \ge N/2} (\mu_k - \mu_j), $$
where $s \le (N-1)/2 - 2$.
Since $ \prod_{l \ne k}
|l - k| = k! (N-k)!$ and $\mu_k - \mu_j = 2(j - k), $
we have that
$$ \prod_{j \ge N/2} |\mu_k - \mu_j | = 2^{(N+1)/2} \frac{(N - k)!}
{( (N-1)/2 - k )! } $$
Define $r = (N-1)/2$.
We need to show the vanishing of
$$\sum_{k = 0}^r (N-2k)^s \frac{(-1)^k (N-k)!}{k! (N-k)! (r - k)! }
$$ for $s \le r - 2$.
It then suffices to show the vanishing of
$$\sum_{k = 0}^r k^s (-1)^k \colvec{r}{k} , \bla \mbox{
$s \le r - 2$.}$$
This follows since one may expand $(1- e^\l)^r$ as a power series
in $e^\l$ using the binomial theorem: we have
$$ (1 - e^\l)^r = \sum_{j= 0}^r (-1)^j \colvec{r}{j}
\exp j \l \: = \sum_{s \ge 0}  \l^s/s!
\sum_{j= 0}^r (-1)^j j^s \colvec{r}{j}
,  $$
but the terms in this expansion corresponding to $s < r$ must vanish
since $1 - e^\l  = \l h(\l)$ for some  function
$h$ of $\l$ which is analytic at $\l = 0$.

\end{document}